\documentclass[preprint,3p,10pt]{elsarticle}

\usepackage{graphicx}			% graphicx package for more complicated commands
\usepackage{amssymb} 			% The amssymb package provides various useful mathematical symbols
\usepackage{amsthm}                       % The amsthm package provides extended theorem environments
\usepackage{mathtools}                    % the mathtools package provides several mathematical functionalities, for example dcases
\usepackage{enumitem}
\usepackage{mathrsfs}
\usepackage[justification=centering]{caption}
\usepackage{subcaption}
\usepackage{adjustbox}
\usepackage{multirow}
\usepackage{mwe}
\usepackage{rotating}
\usepackage{pdflscape}

\biboptions{square,comma,sort&compress}

\newcommand{\ederiv}{\mathrm{d}}

\newcommand{\gradient}{\nabla}		% gradient operator
\newcommand{\curl}{\nabla\times}             % curl operator
\newcommand{\divergence}{\nabla\cdot}   % divergence operator          
\newcommand{\ddt}[2][]{\frac{\ederiv^{#1} #2}{\ederiv t^{#1}}}     % temporal derivative
\newcommand{\pddt}[2][]{\frac{\partial^{#1} #2}{\partial t^{#1}}}   % partial temporal derivative

\newcommand{\hdiv}[1]{H \left( \mathrm{div}, #1 \right) }  % H(div) function space
\newcommand{\hcurl}[1]{H \left( \mathrm{curl}, #1 \right) }  % H(curl) function space
\newcommand{\hgrad}[1]{H \left( \mathrm{grad}, #1 \right) }  % H(grad) function space
\newcommand{\figref}[1]{Figure~\ref{#1}}
\newcommand{\velocity}{\vec{u}}
\newcommand{\velocityd}{\vec{u}_{h}}
\newcommand{\velocitydv}{\boldsymbol{u}}
\newcommand{\velocitydc}[1]{u_{#1}}
\newcommand{\velocitydnext}{\velocityd^{\,k+1}}
\newcommand{\velocitydprevious}{\velocityd^{\,k}}
\newcommand{\velocitydmida}{\tilde{\vec{u}}_{h}^{\,k+\frac{1}{2}}}
\newcommand{\velocitydmidb}{\vec{u}_{h}^{\,k+\frac{1}{2}}}

\newcommand{\velocitydvmidb}{\boldsymbol{u}^{\,k+\frac{1}{2}}}
\newcommand{\velocitydvmidnext}{\boldsymbol{u}_{h}^{\,k+\frac{3}{2}}}
\newcommand{\velocitydmidstart}{\vec{u}_{h}^{\,\frac{1}{2}}}
\newcommand{\velocitydstart}{\vec{u}_{h}^{\,0}}

\newcommand{\vorticity}{\omega}
\newcommand{\vorticityd}{\omega_{h}}
\newcommand{\vorticitydv}{\boldsymbol{\omega}}
\newcommand{\vorticitydc}[1]{\omega_{#1}}
\newcommand{\vorticitydnext}{\vorticityd^{\,k+1}}
\newcommand{\vorticitydprevious}{\vorticityd^{\,k}}
\newcommand{\vorticitydvnext}{\vorticitydv^{\,k+1}}
\newcommand{\vorticitydvprevious}{\vorticitydv^{\,k}}
\newcommand{\vorticitydmida}{\tilde{\omega}_{h}^{\,k+\frac{1}{2}}}

\newcommand{\vorticitydstart}{\vorticityd^{\,0}}

\newcommand{\ufunctionc}{\vec{v}}
\newcommand{\pfunctionc}{q}
\newcommand{\wfunctionc}{\xi}
\newcommand{\uspace}{H(\mathrm{div},\Omega)}
\newcommand{\udspace}{U_{h}}
\newcommand{\ubasis}[1]{\vec{\epsilon}^{\,U}_{#1}}
\newcommand{\udim}{d_{U}}
\newcommand{\wspace}{H(\mathrm{curl},\Omega)}
\newcommand{\wdspace}{W_{h}}
\newcommand{\wbasis}[1]{\epsilon^{W}_{#1}}
\newcommand{\wdim}{d_{W}}
\newcommand{\pspace}{L^{2}(\Omega)}
\newcommand{\pdspace}{Q_{h}}
\newcommand{\pbasis}[1]{\epsilon^{Q}_{#1}}
\newcommand{\pdim}{d_{Q}}
\newcommand{\totalp}{\bar{p}}
\newcommand{\totalpd}{\bar{p}_{h}}
\newcommand{\totalpdv}{\boldsymbol{\bar{p}}}
\newcommand{\totalpdc}[1]{p_{#1}}

\newcommand{\totalpdprevious}{\totalpd^{\,k}}
\newcommand{\totalpdvnext}{\totalpdv^{\,k+1}}

\newcommand{\rtspace}{\mathrm{RT}_{N}}
\newcommand{\qspace}{\mathrm{DG}_{N-1}}
\newcommand{\cgspace}{\mathrm{CG}_{N}}

\newcommand{\matrixoperator}[1]{\boldsymbol{\mathsf{#1}}}
\newcommand{\matrixoperatorc}[1]{\mathsf{#1}}

\journal{Internal note} % the journal name

\begin{document}

	\begin{frontmatter}
		%% Title

         	\title{A conservative, physically compatible, dual-field discretization for turbidity currents: with application to the lock-exchange problem}

         	%% Author and addresses

		\author[tudelft,oxford]{Gonzalo G. de Diego}
		\author[tudelft]{Artur Palha\corref{correspondingAuthor}}
		\author[tudelft]{Marc Gerritsma}
		\address[tudelft]{Delft University of Technology, Faculty of Aerospace Engineering, The Netherlands}
		\address[oxford]{Oxford University, Mathematical Institute, United Kingdom}
		\cortext[correspondingAuthor]{Corresponding author.}
		
		\begin{abstract}
            		In this work we present a structure preserving discretization for turbidity currents based on a mass-, \mbox{energy-,} enstrophy-, and vorticity-conserving formulation for 2D incompressible flows, \cite{palha2016}. This discretization exploits a dual-field formulation for the time evolution of the velocity and vorticity fields together with the transport equation for the particles. Due to its staggered time integration the resulting system of equations is quasi-linear, eliminating the need to solve for a fully nonlinear system of equations. It is shown that this discretization preserves the energy balance equation up to a bounded residual due to the staggering in time of the velocity and vorticity. This leads to a numerical scheme that does not introduce artificial energy dissipation. A comparison with literature results is presented showing that this approach can retrieve the dynamics of the system with a much smaller number of degrees of freedom.
        		\end{abstract}

        		\begin{keyword}
            		Structure preserving \sep Turbidity currents \sep Mixed Finite elements \sep Energy conserving
        		\end{keyword}

	\end{frontmatter}

%% Main text ---------------------------------------------------------------------------------------------------

	\section{Introduction}
		\subsection{Turbidity currents: introduction, relevance, and modelling aproaches}
    			Gravity currents are a general class of phenomena that are recurrent in science and engineering. Essentially, they are horizontally dominated flows driven by hydrostatic pressure gradients. This can be generically pictured as a fluid flowing within another fluid driven by the density difference between the two fluids. Specific cases of gravity currents are present in many situations, e.g.: tunnel fires, CO$_{2}$ sequestration in depleted oil reservoirs, snow avalanches, flow of molten steel on a horizontal surface, lava flows, haboobs (type of sandstorms), doorway flows, and turbidity currents (the focus of this work). For a detailed discussion of gravity currents see for example \cite{Simpson1982, Ungarish2009}, and \cite{Meiburg2015, Huppert2006299} for a more summarized introduction.
    		
    			In the presence of suspended particles, it is frequent to have spatial variations in the bulk density of the fluid. In certain conditions, these variations of the bulk density generate pressure gradients, which induce gravity currents. Turbidity currents are a specific type of these particle-laden gravity currents where the interstitial fluid is a liquid (in most natural cases on Earth it is water) and the particles are typically small (e.g.: sand, clay, etc.) \cite{Meiburg2010, Nasr-Azadani2014}. They occur naturally in lakes and in oceans where they represent a fundamental process for sediment transport \cite{Kneller1999}, since they contribute to both erosion and sedimentation \cite{Meiburg2010}. The range of sediment transport by turbidity currents varies from a few hundreds of meters to thousands of kilometers \cite{Meiburg2010, Nasr-Azadani2014}. Besides their clear geophysical interest, turbidity currents are very relevant from an engineering point of view. First, because they are known to be an important contributor to the formation of deep sea oil reservoirs \cite{Syvitski1996}. Second, due to their significant destructive power, they are one of the main hazards to submarine telecommunication cables and pipelines \cite{parkinson2014, Meiburg2015}.
		
		The study of turbidity currents is particularly challenging. \emph{In situ} investigations are difficult since these particle-laden gravity currents are essentially unpredictable and their shear power typically destroys the measurement devices deployed to analyse them. Small-scale laboratory experiments are a good alternative that provides much understanding into the dynamics of turbidity currents. Unfortunately, they have limitations due to scaling constraints and existing measurement techniques \cite{Kneller2000, parkinson2014}. Another approach, used in the study of gravity currents in general and turbidity currents in particular, is modelling. Given the inherent multitude of flow regimes present in turbidity currents and depending on the particular aspects of interest in the flow, different modelling approaches may be used. These modelling approaches can be grouped into three fundamental types: (i) conceptual models, (ii) depth-averaged models, and (iii) depth-resolving models.
		
		\emph{Conceptual models} for the study of turbidity currents are mainly focussed on establishing analytical estimations (some very accurate) for the front velocity as a function of its height and excess density. These models date back to the seminal work of von K\'{a}rm\'{a}n, \cite{vonKarman1940}, which established the following relation for the Froude number, $F_{h}$,
		\begin{equation}
			F_{h} = \frac{U}{\sqrt{g' h}} = \sqrt{\frac{2}{\sigma}}\,,
		\end{equation}  
		where $U$ is the front velocity of the turbidity current, $h$ represents its height, $g' := \frac{g(\rho_{1} - \rho_{2})}{\rho_{1}}$ stands for the reduced gravity, $\sigma := \frac{\rho_{2}}{\rho_{1}}$ corresponds to the density ratio, and as usual $\rho_{i}$ is the density of the fluid $i$. Over the years, these models have become considerably more advanced and capable of addressing more general cases. For example, in the recent work by Konopliv et al., \cite{Konopliv2016}, a vorticity-based approach is extended to the case of non-Boussinesq gravity currents with success. Although very useful, these models give limited information into the flow and internal processes.
		
		\emph{Depth-averaged models} are a substantial step with regards to the provided insight into the internal dynamics of the flow of gravity currents. The most simple depth-averaged model available is the so called box model. In this approach the gravity current is modelled by a box (a rectangle in 2D, or a rectangular cuboid or cylinder in 3D) that can change its aspect ratio (stretch) as the flow evolves. Although a very simplified model, it can provide useful information on the dynamics of the front and of the height of the gravity current, especially how they are influenced by other parameters of the problem. For a detailed discussion of the box model see for example the monograph by Ungarish \cite{Ungarish2009}. A more detailed depth-averaged representation of gravity currents can be achieved with shallow waters models. When the horizontal dimensions are much larger than the vertical one the shallow waters equations (depth-integrated Navier-Stokes equations) are a good approximation, see for example \cite{Kundu2012} for a derivation. Since the shallow waters equations provide a spatially dependent solution in the horizontal domain, it provides local values of velocity and height of the gravity current. Nevertheless, due to the strong assumptions on the invariance of the flow in the vertical direction this approximation has a limited scope of application and must be used with care in each specific case, as pointed out by Meiburg et al. \cite{Meiburg2015}. For a more detailed discussion and further references see the works by Ungarish \cite{Ungarish2009} and Meiburg et al. \cite{Meiburg2015}.
		
		\emph{Depth-resolved models} are the class of models that have the highest potential of providing accurate and detailed insights into the internal dynamics of the gravity currents, but they are also the most computionally expensive ones. These models solve the full Navier-Stokes equations describing the evolution of the gravity current. For this reason, the full three dimensional information can be retrieved from these models, with a level of detail dictated by the numerical approach used and the available computational resources. Depending on the specific case, different approximations can be used. For example, a Boussinesq approximation is relevant when the driving forces for the gravity current are small density variations; either an Eulerian-Lagrangian formulation tracking each individual particle or a fully Eulerian formulation may be used depending on the particles concentration and dimension. Given the typically high Reynolds numbers present in turbidity currents, a choice must be made regarding the level of detail with which turbulence is modelled. The three main options are (in increasing accuracy and computational cost): Reynolds Averaged Navier-Stokes (RANS), Large-Eddy Simulation (LES), and Direct Numerical Simulation (DNS). RANS equations are obtained by averaging the Navier-Stokes equations in time. Naturally, with this approach the output of these equations are the (temporal) mean fields. RANS simulations rely on a turbulence model that is highly dependent on problem specific parameters. For more details on RANS see Spalart \cite{Spalart2000} for a general introduction, and Meiburg et al. \cite{Meiburg2015} for a summarised discussion focussing on gravity and turbidity currents. The LES formulation is more general and more accurate than RANS, but computationally more expensive. Contrary to RANS, LES relies on a grid dependent spatial filtering: all eddies up to a cutoff scale are resolved, and below this cutoff scale a subgrid-scale model is used to approximate the large scale effects of the small (unresolved) eddies. LES computations have the advantage of producing detailed time-dependent dynamics. For a detailed discussion of LES see Lesieur et al. \cite{Lesieur1996} and Zhiyin \cite{Zhiyin2015}, and Meiburg et al. \cite{Meiburg2015} for a highlight on the aspects related to the simulation of turbidity currents. DNS directly solves the Navier-Stokes equations fully resolving all scales (up to the dissipative Kolmogorov scale). Naturally, this requires very fine grids, which considerably limits the ability to model turbidity currents with very large Reynolds numbers. Nevertheless, the most detailed simulation can only be obtained by DNS simulations.
		
       Early direct numerical simulations of gravity currents focused on the formation of intrusion fronts in lock exchange flows \cite{hartel1997,hartel2000}. The models used in these computations consist of the Boussinesq equation, in which a transport equation models the advection and diffusion of the density variable. Extensions to particle-laden flows were first carried out by Necker et al. \cite{necker2002,necker2005} with the equilibrium Eulerian approach. In these initial computations, the particle velocity field was assumed to be the sum of the carrier fluid's velocity plus a constant settling velocity. Further studies have investigated gravity currents in more complex geometries \cite{blanchette2005,cantero2006} and have considered non-Boussinesq flows \cite{birman2005}. Cantero et al. have extended the investigations on particle-laden flows and accounted for inertial effects, in addition to a settling velocity \cite{cantero2008}. More recently, Parkinson et al. \cite{parkinson2014} have computed particle-laden flows with discontinuous Galerkin finite elements. Espath et al. \cite{espath2014} explore the validity of 2D computations of turbidity currents by comparing them with 3D computations and with experimental data. Meiburg et al. \cite{Meiburg2015} make a thorough review of existing modelling approaches for gravity and turbidity currents and present very detailed 3D simulations, including interaction with obstacles.
       
       In this work, we propose an energy conserving numerical discretization for turbidity currents based on an extension of the previously introduced Mass, Energy, Enstrophy, and Vorticity Conserving (MEEVC) scheme \cite{palha2016, DeDiego2019a}. Due to its conservation properties this method is able to accurately represent the internal energy exchange and to accurately represent the evolution of fluid problems with a smaller number of degrees of freedom. Given its complexity, turbidity currents are a challenging test case for this numerical method. A more efficient numerical discretization can enable the simulation of larger and more complex flows. For the dilute-particle phase a simple approach is implemented: the equilibrium Eulerian method \cite{ferry2001}. This new numerical scheme is assessed by simulating the lock exchange flow. Comparisons are made with existing results in order to evaluate the accuracy of the modified MEEVC scheme.
        
        \subsection{Overview of this work}
			In Section \ref{S2:mc}, a series of modelling possibilities for dilute particle-laden flows are presented, with special emphasis on the equilibrium Eulerian approach because it is the one used in the work. Section \ref{S2:lock_exchange} presents the lock exchange problem, specifying the required boundary conditions for modelling the settling of particles in turbidity currents. Next, in Section \ref{S2:fe}, the finite element discretization of the equilibrium Eulerian approach is presented. The energy budget of a turbidity current is examined in Section \ref{S2:eb} in order to discuss the conservation properties of the solver. Finally, in Section \ref{S2:nr}, numerical results are presented for the lock exchange case and comparisons are made with literature results.

	\section{Modeling considerations}\label{S2:mc}
		In this section, a series of modeling possibilities are presented for dilute suspensions. Of these, the equilibrium Eulerian model is chosen for the construction of the turbidity current solver and the equations of motion are derived. 

		\subsection{Modeling approaches for dilute suspensions}
			A large set of modeling possibilities exist for particle-laden flows. The range of validity of these models is generally determined by two parameters: the Stokes number and the volumetric concentration of particles in the flow \cite{elghobashi1994,balachandar2009}. The Stokes number is defined as the ratio $\tau_p / \tau_f$, where $\tau_p$ and $\tau_f$ are the characteristic time-scale of the particles and the smallest time-scale of the flow (i.e. the Kolmogorov time-scale), respectively. Assuming Stokes flow around the particles, $\tau_p$ can be calculated with,
\begin{equation}\label{eq:tau_p_Stokes}
\tau_p = \frac{\rho_pd^ 2}{18\rho_f\nu},
\end{equation}

		\noindent where $\rho_p$ and $\rho_f$ are the particle and fluid densities, respectively, $d$ is the particle diameter and $\nu$ the dynamic viscosity of the fluid \cite{elghobashi1994}.

\begin{figure}[]
			\centering
				\includegraphics[width=0.45\textwidth, trim=0.0 100.0 0.0 100.0, keepaspectratio]{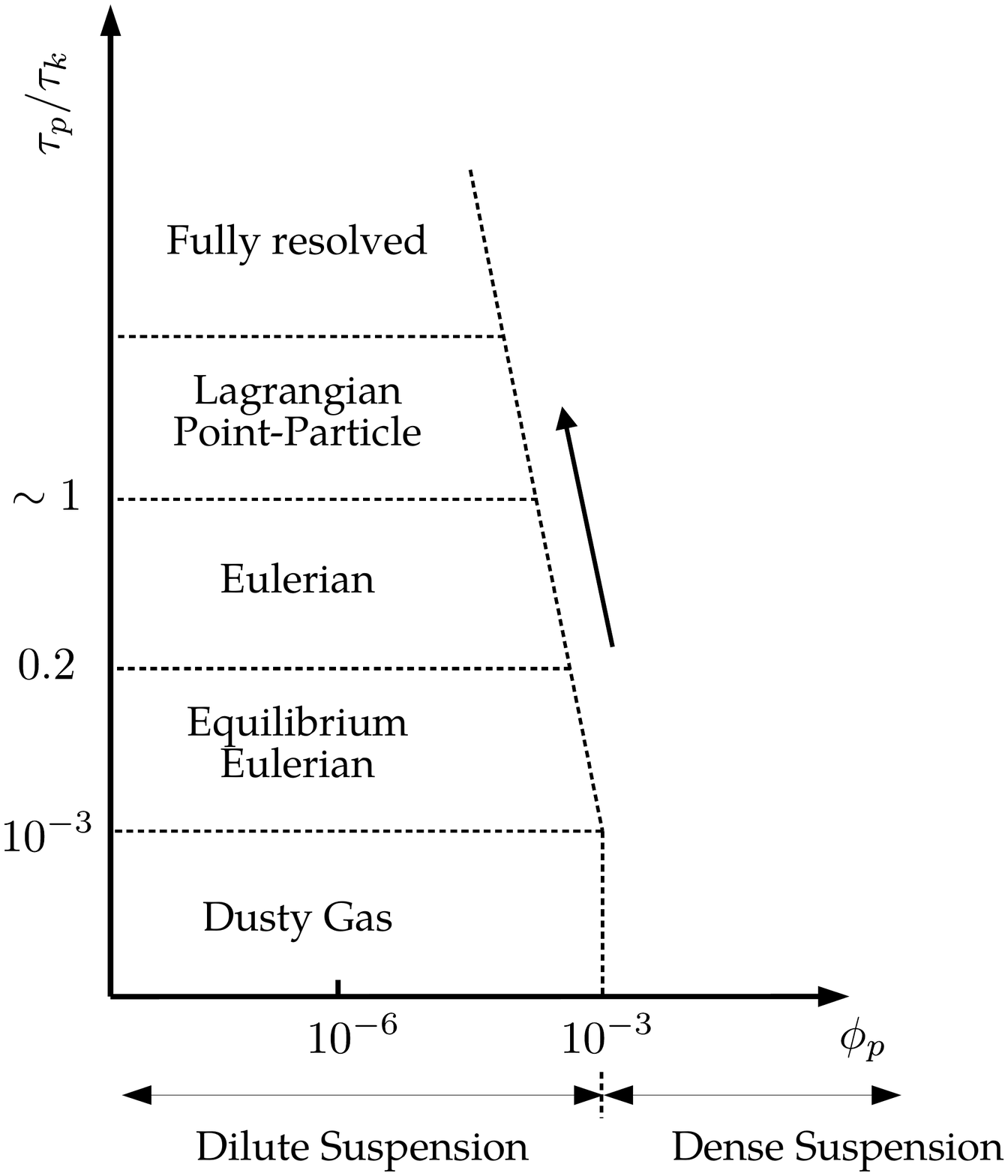}
				\caption{Different approaches to modeling particle-laden flows according to the volume fraction $\phi_p$ of the particle phase and time scale ratio (Stokes number). The direction indicated by the vertical arrow reflects the increasing importance of particle collisions for larger particle sizes. Figure taken from Balachandar \cite{balachandar2009}.}
				\label{fig:particle_modeling}
		\end{figure}

		Figure \ref{fig:particle_modeling} summarizes the modeling approaches for dilute particle suspensions. In this case, the volumetric concentration of particles is limited to $\phi_p < 10^{-3}$, such that a two-way coupling based on momentum exchange exists between the fluid and the particles. For higher concentrations, the flow is considered a dense suspension and particle-particle interaction becomes important (four-way coupling), \cite{elghobashi1994}.

		The different modeling approaches indicated in Figure \ref{fig:particle_modeling} are briefly explained below.
		\begin{itemize}
			\item \textbf{Dusty gas approach}. Particles are assumed to be sufficiently small such that they follow the carrier fluid. The equations to be solved are the Navier-Stokes equations with a modified density together with a transport equation for particles.

			\item \textbf{Equilibrium Eulerian approach}. Developed by Ferry and Balachandar and also known as the \emph{fast Eulerian method} \cite{ferry2001}, this method retains the simplicity of the dusty gas approach by assuming that the particle velocity field can be expanded in terms of the carrier fluid velocity field and calculated by means of an algebraic equation. This method gives adequate results for particles with Stokes numbers just below unity and density ratio parameters $\rho_p / \rho_f$ of order $\mathcal{O}(1)$ or lower. Further attention is given to this method in Section \ref{ss2:eqEul}.

			\item \textbf{Eulerian approach}. Both the fluid and the particles are treated as two different continuum fluid phases. In this case, momentum and energy equations must also be solved for the particle phase. The principal restriction for this method is that within a small volume of fluid the particles must have a unique velocity, momentum and energy.

			\item \textbf{Lagrangian approach}. For $\tau_p/\tau_f > 1$, the relative particle sizes are such that the continuum phase assumption is no longer realistic, that is, assuming uniqueness of particle properties over small areas of the fluid is not valid. In the Lagrangian approach, particles are represented in a Lagrangian reference frame and their properties are calculated with probabilistic methods. In this case, there is no requirement for uniqueness and the flow equations are valid until the Stokes number is such that the point-particle assumption no longer makes sense.
		\end{itemize}

		In this work, the equilibrium Eulerian approach will be used due to its simplicity and its capacity to capture both inertial and gravitational effects characteristic of turbidity currents. In the following subsection, the equilibrium assumption will be applied to the Eulerian-Eulerian equations and a set of equations for dilute suspensions will be derived.

		\subsection{Equilibrium Eulerian equations for a dilute suspension}\label{ss2:eqEul}

			A set of equations for a fluid-particle system studied from an Eulerian point of view can be obtained by mixture theory or by ensemble averaging \cite{joseph1990}. The volumetric concentrations of the fluid and particle phases are denoted $\phi_f$ and $\phi_p$ respectively, the densities $\rho_f$ and $\rho_p$ and the velocities $\vec{u}_p$ and $\vec{u}_f$. 

			By assuming constant densities for both phases and no mass transfer, the equations of motion in dimensional form can be written as follows \cite{cantero2008},
			\begin{equation}\label{eq:euler_fluid_particle}
				\begin{dcases}
					\frac{\partial \phi_f}{\partial t} + \divergence \left( \phi_f \vec{u}_f \right) = 0, \\
					\frac{\partial \phi_p}{\partial t} + \divergence \left( \phi_p \vec{u}_p \right) = 0, \\
					\phi_p\rho_p\frac{\mathrm{D}_p \vec{u}_p}{\mathrm{D} t}= \phi_p\left(\rho_p - \rho_f\right) \vec{g} - \phi_d\gradient p + \vec{F},\\
					\phi_f\rho_f\frac{\mathrm{D}_f \vec{u}_f}{\mathrm{D} t}=-\phi_f\gradient p + \mu\Delta\vec{u}_v - \vec{F},
				\end{dcases}
			\end{equation}

			\noindent where $\vec{F}$ represents the net hydrodynamic interaction between both phases, $\vec{g} = g \vec{e}_g$ is the gravitational acceleration vector ($\vec{e}_g$ is the unit vector pointing in the direction of the gravitational force) and $\vec{u}_v = \phi_f \vec{u}_f + \phi_p \vec{u}_p$ is the composite velocity. $\mathrm{D}_f / \mathrm{D} t$ and $\mathrm{D}_p / \mathrm{D} t$ refer to the material derivatives following the fluid and particle velocity fields, respectively. The volume fractions satisfy the identity $\phi_f + \phi_p = 1$; hence, by adding the first two equations of \eqref{eq:euler_fluid_particle},
			\begin{equation}\label{eq:uv_div}
				\divergence \vec{u}_v = \divergence \left( \phi_f \vec{u}_f + \phi_p \vec{u}_p \right) = 0.
			\end{equation}

			In the equilibrium Eulerian model, the velocity field of the particle phase, $\vec{u}_p$, is expanded about $\tau_p = 0$ and defined in terms of $u_f$. Given the density ratio parameter $\beta = 3/(2\rho_p/\rho_f + 1)$, for $\beta \sim \mathcal{O}(1)$, the following approximation holds up to $\mathcal{O}(\tau_p^{3/2})$:
\begin{equation}\label{eq:p_vel_1}
\vec{u}_p \approx \vec{u}_f + \tau_p (1-\beta) \left( \vec{g} - \frac{\mathrm{D}_f \vec{u}_f}{\mathrm{D} t} \right),
\end{equation}

			For simplicity and to conform with existing computations \cite{necker2002,necker2005,parkinson2014,espath2014} for validation purposes, the inertial term in \eqref{eq:p_vel_1} is neglected, so that,
			\begin{equation}\label{eq:p_vel_2}
				\vec{u}_p \approx \vec{u}_f + \tau_p (1-\beta) \vec{g} = \vec{u}_f + u_s \vec{e}_g,
			\end{equation}

			\noindent where $u_s = \tau_p (1-\beta) g$ is the settling velocity of the particles.

			The fact that the particle velocity field is completely defined by the velocity field of the carrier fluid leads to the following question: how will the initial particle distribution affect the resulting particle velocity field? It is possible to prove that, for sufficiently small particle sizes, an \emph{equilibrium} particle field exists such that it depends only on fluid quantities \cite{ferry2001}. In essence, this means that the transients arising from the initial conditions decay exponentially fast when the characteristic time-scale $\tau_p$ is much smaller than that of the fluid $\tau_f$.

			It is important to take into account that certain circumstances exist in which the Stokes number $\tau_p/\tau_f$ is small but the equilibrium condition does not hold. An example is the injection of particles into a fluid domain; near the point of injection, particle motion will be dominated by the injection process and only after an initial transient will it reach equilibrium with respect to the carrier fluid \cite{balachandar2009}. Furthermore, for increasing particle concentrations, the mean time between particle collisions will decrease until it reaches the particle response time $\tau_p$ and exerts a considerable effect on the particle velocity field.

			Implicit in \eqref{eq:p_vel_2} is the assumption:
			\begin{equation}\label{eq:Dp_1}
				\frac{\mathrm{D}_p\vec{u}_p}{\mathrm{D} t} \approx \frac{\mathrm{D}_f\vec{u}_f}{\mathrm{D} t}.
			\end{equation}

			Taking into account \eqref{eq:Dp_1}, the momentum equations in \eqref{eq:euler_fluid_particle} can be combined in order to obtain
			\begin{equation}\label{eq:momentum_fluid_part}
				\left( \phi_f\rho_f + \phi_p\rho_p \right) \frac{\mathrm{D}_f \vec{u}_f}{\mathrm{D} t}=-\gradient p + \mu\Delta\vec{u}_v + \phi_p\left(\rho_p - \rho_f\right) \vec{g}.
			\end{equation}

			Considering that dilute suspensions are to be modeled and therefore $\phi_p << 1$, the composite velocity can be considered equal to the fluid velocity field, that is, $\vec{u}_v \approx \vec{u}_f$. By using \eqref{eq:uv_div}, the conservation of mass for the fluid phase can be described with a divergence-free constraint, 
			\begin{equation}\label{eq:u_div}
				\divergence \vec{u}_f \approx 0.
			\end{equation}

			Finally, the Boussinesq approximation is assumed to hold, such that density variations are considered to be small and only influence the bouyancy term of \eqref{eq:momentum_fluid_part}. The density term on the left-hand side of \eqref{eq:momentum_fluid_part} is assumed constant, such that $\left( \phi_f\rho_f + \phi_p\rho_p \right) \approx \rho_f$. In order to write the simplified equations of motion in non-dimensional form, a characteristic velocity, length and density are defined. The characteristic density is taken as $\rho_f$ and, when considering a lock exchange flow, the characteristic length is taken as the height of the channel, $H$, see Figure \ref{fig:problem}. The buoyancy velocity, defined as
			\begin{equation}\label{eq:ub}
				u_b = \sqrt{gH\frac{\rho_p-\rho_f}{\rho_f}\phi_{p,\mathrm{max}}},
			\end{equation}

			\noindent is used as the characteristic velocity. The term $\phi_{p,\mathrm{max}}$ represents the maximum volumetric concentration of particles in the domain at the initial time instant. For simplicity of notation, the non-dimensional fluid velocity, $\vec{u}_f$, will be represented with $\vec{u}$ and the variable $\phi = \phi_p / \phi_{p,\mathrm{max}}$ will be used as a normalized measure of particle concentration. The non-dimensional equations of motion are given by,
			\begin{equation}\label{eq:NS_particle}
				\begin{dcases}
					\divergence \vec{u} = 0,\\
					\pddt{\vec{u}} + \left( \vec{u} \cdot \nabla \right)\ \vec{u} =-\gradient p + \frac{1}{\sqrt{\mathrm{Gr}}} \Delta\vec{u} + \phi \vec{e}_g,\\
					\pddt{\phi} + \left( \vec{u} + u_s\vec{e}_g  \right) \cdot \gradient \phi = \frac{1}{\sqrt{\mathrm{Gr} \mathrm{Sc}^2}} \Delta \phi.
				\end{dcases}
			\end{equation}

			In Equation \eqref{eq:NS_particle}, a diffusion term has been added to the equation for conservation of particles. This is common practice when considering particle-laden flows, as it avoids the formation of sharp concentrations of particles that could lead to numerical instabilities \cite{necker2002,necker2005,cantero2008}. Furthermore, certain authors indicate that this term has a physical significance, such as the spreading of particles in time due to hydrodynamic diffusion \cite{necker2005} or the departure of the equilibrium assumption due to close interaction of particles \cite{cantero2008}.

			The two non-dimensional parameters that appear in \eqref{eq:NS_particle} are the Grashof number $\mathrm{Gr}$ and the Schmidt number $\mathrm{Sc}$, defined by
\begin{equation}\label{eq:Gr}
\mathrm{Gr} = \left( \frac{u_bH}{\nu} \right)^2 \quad \text{and} \quad \mathrm{Sc} = \frac{\nu}{\kappa}.
\end{equation}

			The Grashof number represents the ratio of buoyancy forces and viscous forces and it is proportional to $\mathrm{Re}^2$. The Schmidt number compares the viscous diffusivity of the fluid to the molecular diffusivity of the particle field, given by $\kappa$.
	
	\section{The 2D lock exchange flow}\label{S2:lock_exchange}
		In this section, the 2D lock-exchange flow problem of a particle-driven current is presented together with the boundary conditions that will be prescribed in order to capture the deposition of sediment along the lower boundary. A lock-exchange flow consists of the mutual propagation of two flows with different densities which were initially separated by a membrane. Many physically relevant features of gravity currents can be observed in this type of flow, such as the body-head structure that gradually develops as the flow evolves \cite{hartel1997,hartel2000,cantero2008}. Extensive experimental studies have been carried out with lock-exchange flows \cite{simpson1997}; more recently, direct numerical simulations have been used in order to capture more subtle physical aspects of the flow or to test numerical solvers for gravity currents \cite{hartel1997,hartel2000,necker2002,birman2005,blanchette2005,cantero2006,cantero2008,espath2014,parkinson2014}.

\begin{figure}[]
			\centering
				\includegraphics[width=0.7\textwidth, keepaspectratio]{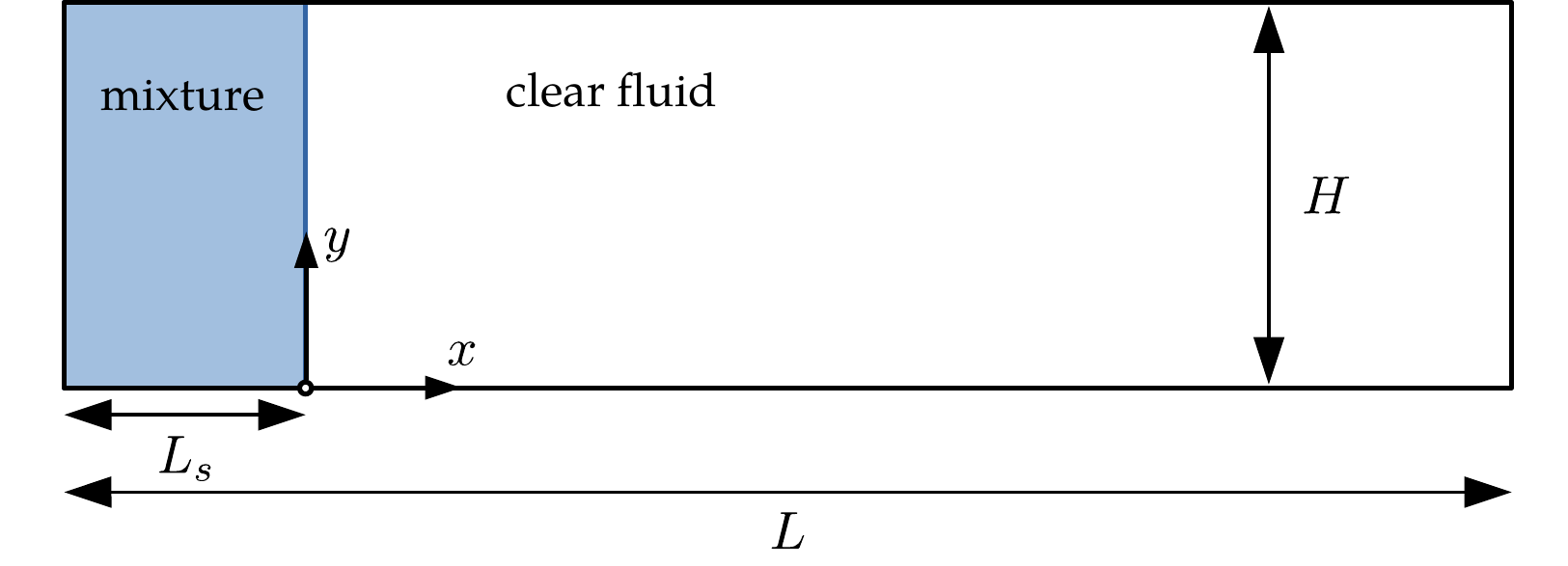}
				\caption{Initial set-up of the lock-exchange test-case to be computed in this work.}
				\label{fig:problem}
		\end{figure}

		The initial set-up of the lock-exchange flow computed in this work is shown in Figure \ref{fig:problem} and corresponds to the one presented by Necker et al. \cite{necker2002} and used in the work of Parkinson et al. \cite{parkinson2014} and Espath et al. \cite{espath2014}. The domain consists of a plane channel of height $H$ and length $L$. At $t=0$, a mixture of the particle phase and the fluid phase is contained in the region $\left[ -L_s, 0 \right] \times \left[ 0, H \right]$. For $t > 0$, the mixture is put in touch with the clear fluid and an intrusion front develops as the heavier fluid propagates along the lower boundary.

		\subsection{Boundary conditions}\label{subsec:le_bc}
\begin{figure}[]
				\centering
					\includegraphics[width=0.7\textwidth, keepaspectratio]{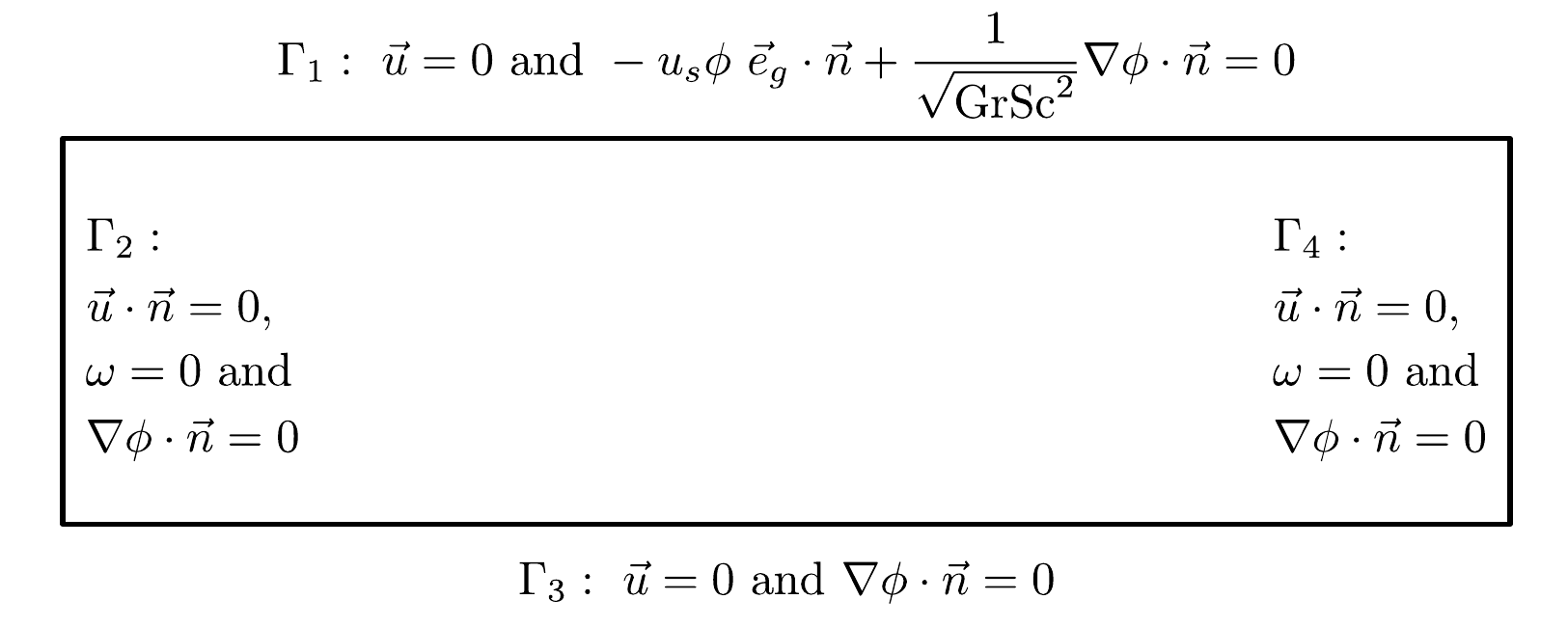}
					\caption{Boundary conditions for the lock exchange test-case.}
					\label{fig:bc}
			\end{figure}

			In order to solve \eqref{eq:NS_particle}, boundary conditions must be given for the velocity $\vec{u}$ and particle concentration $\phi$. Boundary conditions will be prescribed in the same way as Espath et al. \cite{espath2014}, such that the lower wall allows for the convection of the flow particles at a settling velocity $u_s$. To this end, the boundary $\partial \Omega$ is partitioned into 4 regions which correspond with the 4 faces of the rectangle, such that $\partial \Omega =  \cup\ \Gamma_i$, for $i = 1,2,3$ and $4$, see Figure \ref{fig:bc}.

			The velocity boundary conditions consist of no-slip and no-penetration conditions along the upper and lower walls, such that 
			\begin{equation}\label{eq:noslip_bc}
				\vec{u} = 0 \quad \text{on}\ \Gamma_1 \cup \Gamma_3,
			\end{equation}

			\noindent and shear free conditions with no-penetration along the lateral walls,
			\begin{equation}\label{eq:slip_bc}
				\vec{u} \cdot \vec{n} = 0\quad \text{and}\quad \gradient \vec{u} \cdot \vec{n} = 0 \quad \text{on}\quad \Gamma_2 \cup \Gamma_4.
			\end{equation}

			Taking into account that \eqref{eq:NS_particle} will be solved in a $(\vec{u}, \omega, p)$ formulation, the shear-free boundary conditions can also be written in the following form,
			\begin{equation}\label{eq:slip_bc2}
				\vec{u} \cdot \vec{n} = 0 \quad \text{and} \quad \omega = 0 \quad \text{on}\quad \Gamma_2 \cup \Gamma_4.
			\end{equation}

			A fundamental feature of turbidity currents is the deposition and resuspension of particles, which modifies the driving force of the current \cite{cantero2008}. For the lock exchange flow, deposition of particles is assumed to occur on the lower boundary, $\Gamma_3$, due to advection with $u_s \vec{e}_g$. Over $\Gamma_1 \cup \Gamma_2 \cup \Gamma_4$, no particle transport is assumed to hold.

			In order to understand how these boundary conditions can be prescribed, the transport equation for $\phi$ must be integrated over the domain $\Omega$. Taking into account that $\divergence \vec{u} = 0$ in $\Omega$ and $\vec{u} \cdot \vec{n} = 0$ on $\partial \Omega$, integrating the particle transport equation in \eqref{eq:NS_particle} over $\Omega$ yields the following result after applying the divergence theorem,
			\begin{equation}\label{eq:p_flux}
				\ddt{} \int_{\Omega} \phi\ \mathrm{d}\Omega = \int_{\partial \Omega} \left( -\phi u_s \vec{e}_g + \frac{1}{\sqrt{\mathrm{Gr} \mathrm{Sc}^2}} \gradient \phi \right) \cdot \vec{n}\ \mathrm{d}\Gamma.
			\end{equation}

			Therefore, if the boundary condition for $\phi$ on $\Gamma_1 \cup \Gamma_2 \cup \Gamma_4$ is zero particle flux, the integral in Equation \eqref{eq:p_flux} must be equal to zero, such that 
			\begin{equation}\label{eq:p_no_flux_bc}
				\left( -\phi u_s \vec{e}_g + \frac{1}{\sqrt{\mathrm{Gr} \mathrm{Sc}^2}} \gradient \phi \right) \cdot \vec{n} = 0 \quad \text{on}\ \Gamma_1 \cup \Gamma_2 \cup \Gamma_4.
			\end{equation}

			If the gravity vector is set to $\vec{e}_g = (0,-1)$, the following boundary conditions hold:
			\begin{equation}\label{eq:p_no_flux_bc2}
				\begin{dcases}
					- u_s \phi\ \vec{e}_g \cdot \vec{n} + \frac{1}{\sqrt{\mathrm{Gr} \mathrm{Sc}^2}} \gradient \phi \cdot \vec{n} = 0 & \text{on}\ \Gamma_1, \\
					\gradient \phi \cdot \vec{n} = 0&  \text{on}\ \Gamma_2 \cup \Gamma_4.
				\end{dcases}
			\end{equation}

			Deposition of particles is enabled by assuming that the boundary condition over $\Gamma_2$ is a particle flux at a velocity $u_s$. In order to do so, the dissipation term in \eqref{eq:p_flux} is eliminated, such that
			\begin{equation}\label{eq:p_flux_bc}
				\gradient \phi \cdot \vec{n} = 0 \quad  \text{on}\ \Gamma_3.
			\end{equation}

			Therefore, given the boundary conditions \eqref{eq:p_no_flux_bc2} and \eqref{eq:p_flux_bc}, the variation of the total particle concentration is determined with the following expression,
			\begin{equation}\label{eq:p_flux_2}
				\ddt{} \int_{\Omega} \phi\ \mathrm{d}\Omega = - \int_{\Gamma_3} \phi u_s\ \mathrm{d}\Gamma.
			\end{equation}

	\section{Numerical discretization} \label{S2:fe}
		In this work we introduce a structure preserving discretization for turbidity currents based on the MEEVC discretization,  \cite{palha2016, DeDiego2019a}. We start by briefly introducing the MEEVC discretization for a homogeneous fluid first for the case of periodic boundaries and then for the case of non-periodic boundary conditions. The MEEVC method will be presented on an \emph{ad hoc} basis, for a detailed discussion on the construction of the method and proof of its properties the reader is directed to  \cite{palha2016, DeDiego2019a}. Once the MEEVC method for a homogeneous fluid has been presented the extension to turbidity currents will be discussed in detail.
		
		\subsection{Brief introduction to the MEEVC discretization}
			\subsubsection{Spatial discretization (periodic boundary conditions)}
				The most common form of the Navier-Stokes equations for 2D viscous incompressible flows, as seen above, is expressed as a set of conservation laws for momentum and mass involving the velocity $\velocity: \Omega\times (0, t_{F}]\mapsto\mathbb{R}^{2}$ and pressure $p: \Omega\times (0, t_{F}]\mapsto\mathbb{R}$:
		\begin{equation}
			\begin{dcases}
				\frac{\partial\velocity}{\partial t} + \left(\velocity\cdot\nabla\right)\velocity + \nabla p = \nu\Delta\velocity + \vec{s}, & \text{in }\Omega\times (0, t_{F}]\,,\\
				\nabla\cdot\velocity = 0, & \text{in }\Omega\times (0, t_{F}]\,,
			\end{dcases} \label{eq::ns_convective_form}
		\end{equation}
		with the final time instant $t_{F} > 0$, $\nu$ the kinematic viscosity, $\vec{s}$ the body force per unit mass, and $\Delta = \nabla\cdot\nabla$ the Laplace operator. These equations are valid on the fluid domain $\Omega$, together with suitable initial and boundary conditions. At this stage, we consider only periodic boundary conditions and we set $\vec{s}=0$.
		
		The form of the Navier-Stokes equations presented in \eqref{eq::ns_convective_form} is the so called \emph{convective form}. Its name stems from the particular form of the nonlinear term $ \left(\velocity\cdot\nabla\right)\velocity$, which underlines its convective nature. This form is not unique. Using well known vector calculus identities it is possible to rewrite the nonlinear term in three other forms: \emph{divergence form},\emph{ skew-symmetric form}, and \emph{rotational form} (see for example Zang \cite{Zang1991}, Morinishi \cite{Morinishi1998}, and R{\o}nquist \cite{Ronquist1996}).
				
				The MEEVC discretization uses as a starting point the \emph{rotational form}, which makes use of the vorticity $\vorticity:=\nabla\times\velocity$:
		\begin{equation}
			\begin{dcases}
				\frac{\partial\velocity}{\partial t} + \omega\times\velocity + \nabla \totalp = -\nu\nabla\times\omega,  & \text{in }\Omega\times (0, t_{F}]\,,\\
				\nabla\cdot\velocity = 0, & \text{in }\Omega\times (0, t_{F}]\,,
			\end{dcases} \label{eq::ns_rotational_form}
		\end{equation}
		where the static pressure $p$ is replaced by the total pressure $\totalp := \frac{1}{2}\velocity\cdot\velocity + p$.
			
			By introducing an unknown, the vorticity, we lack one additional equation in order to close the system. Different options could be used (each with its advantages and drawbacks): (i) $\omega := \nabla\times\velocity$ (kinematic equation), (ii) $\Delta \psi = \omega$ with $\nabla\times\psi = \velocity$ (stream function), and (iii) a dynamic equation for vorticity. In the MEEVC method the third route is followed. This route, referred to as dual-field formulation, employs two dynamic equations for the computation of the evolution of the flow. As opposed to a traditional approach where the vorticity is related to a statically defined curl operator, in this work this relation is dynamic and evolves in time. This is the key component to obtain the conservation properties present in this numerical disretization. For a detailed discussion the reader is directed to \cite{palha2016}.
			
			It is well known that by taking the curl of the momentum equation and using the kinematic definition $\vorticity := \nabla\times\velocity$ we can obtain the flow equations based on vorticity transport in skew-symmetric form:
		\begin{equation}
			\begin{dcases}
				\frac{\partial\vorticity}{\partial t} + \frac{1}{2}\left(\velocity\cdot\nabla\right)\vorticity + \frac{1}{2}\nabla\cdot\left(\velocity\,\vorticity\right) = \nu\Delta\vorticity,  & \text{in }\Omega\times (0, t_{F}]\,,\\
				\nabla\cdot\velocity = 0, & \text{in }\Omega\times (0, t_{F}]\,,\\
				\vorticity = \nabla\times\velocity\,,  & \text{in }\Omega\times (0, t_{F}]\,.
			\end{dcases} \label{eq::ns_vorticity_transport}
		\end{equation}
		This velocity-vorticity $(\velocity,\vorticity)$ formulation of the Navier-Stokes equations is of particular interest for vortex dominated flows, see for example Gatski \cite{Gatski1991} for an overview and Daube \cite{Daube1992} and Clercx \cite{Clercx1997a} for applications.
		
		The MEEVC discretization starts with a $(\velocity,\vorticity)$ formulation by combining the rotational form \eqref{eq::ns_rotational_form} with the vorticity transport equation \eqref{eq::ns_vorticity_transport}:
		\begin{equation}
			\begin{dcases}
				\frac{\partial\velocity}{\partial t} + \omega\times\velocity + \nabla \totalp = -\nu\nabla\times\vorticity,  & \text{in }\Omega\times (0, t_{F}]\,,\\
				\frac{\partial\vorticity}{\partial t} + \frac{1}{2}\left(\velocity\cdot\nabla\right)\vorticity + \frac{1}{2}\nabla\cdot\left(\velocity\,\vorticity\right) = \nu\Delta\vorticity,  & \text{in }\Omega\times (0, t_{F}]\,,\\
				\nabla\cdot\velocity = 0 & \text{in }\Omega\times (0, t_{F}]\,.
			\end{dcases} \label{eq::ns_meevc_form}
		\end{equation}
	
		An important aspect we wish to stress is that although at the continuous level the kinematic definition $\vorticity := \nabla\times\velocity$ is exactly satisfied, at the discrete level it is not always guaranteed that this identity holds exactly. In fact, in the discretization presented here this identity is satisfied only approximately. This, as can be seen in \cite{palha2016, DeDiego2019a}, enables the construction of a mass, energy, enstrophy and vorticity conserving discretization (MEEVC). 
		
		The next step for developing the MEEVC discretization is the construction of the weak form of \eqref{eq::ns_meevc_form}, as is standard in finite elements:
		\begin{equation}
			\begin{dcases}
				\text{Find } \velocity\in \uspace, \totalp\in \pspace \text{ and } \vorticity \in \wspace \text{ such that:}\\
				\langle\frac{\partial\velocity}{\partial t},\ufunctionc\rangle_{\Omega} + \langle\omega\times\velocity,\ufunctionc\rangle_{\Omega} - \langle \totalp,\nabla\cdot\ufunctionc\rangle_{\Omega} = -\nu\langle\nabla\times\vorticity,\ufunctionc\rangle_{\Omega}, & \forall \ufunctionc\in \uspace, \\
				\langle\frac{\partial\vorticity}{\partial t},\wfunctionc\rangle_{\Omega} - \frac{1}{2}\langle\vorticity,\nabla\cdot\left(\velocity\,\wfunctionc\right)\rangle_{\Omega} +  \frac{1}{2}\langle\nabla\cdot\left(\velocity\,\vorticity,\right),\wfunctionc\rangle_{\Omega} = \nu\langle\nabla\times\vorticity,\nabla\times\wfunctionc\rangle_{\Omega}, & \forall \wfunctionc\in \wspace, \\
				\langle\nabla\cdot\velocity,\pfunctionc\rangle_{\Omega} = 0, &\forall \pfunctionc\in \pspace\,,
			\end{dcases} \label{eq::ns_meevc_weak_form_continuous}
		\end{equation}
		where we have used integration by parts and the periodic boundary conditions to obtain the identities $\langle \totalp,\nabla\cdot\ufunctionc\rangle_{\Omega} = - \langle \nabla\totalp,\ufunctionc\rangle_{\Omega}$,   $\langle\vorticity,\nabla\cdot\left(\velocity\,\wfunctionc\right)\rangle_{\Omega} = -\langle\left(\velocity\cdot\nabla\right)\vorticity,\wfunctionc\rangle_{\Omega}$ and $\langle\nabla\times\vorticity,\nabla\times\wfunctionc\rangle_{\Omega} = \langle\Delta\vorticity,\wfunctionc\rangle_{\Omega}$. The space $\pspace$ corresponds to square integrable functions and the spaces $\uspace$ and $\wspace$ contain square integrable functions whose divergence and curl are also square integrable. The angled brackets $\langle \cdot, \cdot \rangle_\Omega$ denoted the $L^2$ inner product in $\Omega$.
			
			The crucial step to transform these infinite dimensional continuous equations into computable finite dimensional equations relies on a choice of adequate conforming finite dimensional function spaces, where we will seek our discrete solutions for velocity $\velocityd$, pressure $\totalpd$ and vorticity $\vorticityd$:
		\begin{equation}
			\velocityd \in \udspace \subset \uspace, \quad \totalpd \in \pdspace \subset \pspace \quad \mathrm{and} \quad \vorticityd \in \wdspace \subset \wspace.
		\end{equation}
		
		As usual, each of these finite dimensional function spaces, $\udspace$, $\pdspace$ and $\wdspace$, has an associated finite set of basis functions, $\ubasis{i}$, $\pbasis{i}$, $\wbasis{i}$, such that
		\begin{equation}
			\udspace = \mathrm{span}\{\ubasis{1}, \dots,\ubasis{\udim}\}, \quad \pdspace = \mathrm{span}\{\pbasis{1}, \dots,\pbasis{\pdim}\}\quad\mathrm{and}\quad\wdspace = \mathrm{span}\{\wbasis{1}, \dots,\wbasis{\wdim}\},
		\end{equation}
		 where $\udim$, $\pdim$ and $\wdim$ denote the dimension of the discrete function spaces and therefore correspond to the number of degrees of freedom for each of the unknowns. As a consequence, the approximate solutions for velocity, pressure and vorticity can be expressed as a linear combination of these basis functions
		 \begin{equation}
		 	\velocityd := \sum_{i=1}^{\udim}\velocitydc{i}\,\ubasis{i}, \quad \totalpd := \sum_{i=1}^{\pdim}\totalpdc{i}\,\pbasis{i} \quad\mathrm{and}\quad \vorticityd := \sum_{i=1}^{\wdim}\vorticitydc{i}\,\wbasis{i}, \label{eq:basis_expansion}
		 \end{equation}
		 with $\velocitydc{i}$, $\totalpdc{i}$ and $\vorticitydc{i}$ the degrees of freedom of velocity, total pressure and vorticity, respectively. Since the Navier-Stokes equations form a time dependent set of equations, in general these coefficients will be time dependent, $\velocitydc{i} = \velocitydc{i}(t)$, $\totalpdc{i}=\totalpdc{i}(t)$ and $\vorticitydc{i}=\vorticitydc{i}(t)$.
		 
		The choice of the finite dimensional function spaces dictates the properties of the discretization. In order to have exact conservation of mass, energy, enstrophy, and total vorticity we must choose these function spaces such that they form a Hilbert subcomplex
		\begin{equation}
			0 \longrightarrow \wdspace \stackrel{\nabla\times}{\longrightarrow}\udspace \stackrel{\nabla\cdot}{\longrightarrow}\pdspace \longrightarrow 0\,,
		\end{equation}
		that mimics the 2D Hilbert complex associated to the continuous function spaces:
		\begin{equation}
			0 \longrightarrow \wspace \stackrel{\nabla\times}{\longrightarrow}\uspace \stackrel{\nabla\cdot}{\longrightarrow}\pspace \longrightarrow 0\,.
		\end{equation}
		The Hilbert complex is an important structure that is intimately related to the de Rham complex of differential forms. The construction of a discrete subcomplex is an important requirement to obtain stable and accurate finite element discretizations, see for example \cite{arnold2010finite,Palha2014,{bossavit_japan_computational_1,bossavit_japan_computational_2,bossavit_japan_computational_3,bossavit_japan_computational_4,bossavit_japan_computational_5}} for a general discussion and  \cite{palha2016} for the specific discussion on the MEEVC formulation.
			
			One specific choice of discrete function spaces is
			\begin{equation}
				\wdspace = \cgspace, \quad \udspace = \rtspace, \quad\text{and}\quad \pdspace=\qspace\,,
			\end{equation}
			where $\cgspace$ are the Lagrange elements of degree $N$, see \cite{kirby2012}, $\rtspace$ are the Raviart-Thomas elements of degree $N$, see \cite{RaviartThomas1977,kirby2012}, and $\qspace$ are the discontinuous Lagrange elements of degree $(N-1)$, see \cite{kirby2012}.
			
			Using the discrete expansions for $\velocityd$, $\totalpd$ and $\vorticityd$, \eqref{eq:basis_expansion}, can be transformed into its discrete counterpart
		\begin{equation}
			\begin{dcases}
				\text{Find } \velocitydv\in \mathbb{R}^{\udim}, \totalpdv\in \mathbb{R}^{\pdim} \text{ and } \vorticitydv \in \mathbb{R}^{\wdim} \text{ such that:}\\
				\sum_{i=1}^{\udim}\frac{\mathrm{d}\velocitydc{i}}{\mathrm{d}t}\langle\ubasis{i},\ubasis{j}\rangle_{\Omega} + \sum_{i=1}^{\udim}\velocitydc{i}\langle\vorticityd\times\ubasis{i},\ubasis{j}\rangle_{\Omega} - \sum_{k=1}^{\pdim}\totalpdc{k}\langle \pbasis{k},\nabla\cdot\ubasis{j}\rangle_{\Omega} = -\nu\langle\nabla\times\vorticityd,\ubasis{j}\rangle_{\Omega},\quad j=1,\dots,\udim, \\
				\sum_{i=1}^{\wdim}\frac{\mathrm{d}\vorticitydc{i}}{\mathrm{d}t}\langle\wbasis{i},\wbasis{j}\rangle_{\Omega} - \sum_{i=1}^{\wdim}\frac{\vorticitydc{i}}{2}\langle\wbasis{i},\nabla\cdot\left(\velocityd\,\wbasis{j}\right)\rangle_{\Omega} + \sum_{i=1}^{\wdim}\frac{\vorticitydc{i}}{2}\langle\nabla\cdot\left(\velocityd\,\wbasis{i}\right),\wbasis{j}\rangle_{\Omega} = \\
				\qquad\qquad\qquad\qquad\qquad\qquad\qquad\qquad\qquad\qquad\qquad\qquad\nu\sum_{i=1}^{\wdim}\vorticitydc{i}\langle\nabla\times\wbasis{i},\nabla\times\wbasis{j}\rangle_{\Omega}, \quad j=1,\dots,\wdim, \\
				\sum_{i=1}^{\udim}\velocitydc{i}\langle\nabla\cdot\ubasis{i},\pbasis{j}\rangle_{\Omega} = 0, \quad j = 1,\dots,\pdim\,,
			\end{dcases} \label{eq::ns_meevc_weak_form_discrete_expansion}
		\end{equation}
		with $\velocitydv := [\velocitydc{1},\dots,\velocitydc{\udim}]^{\top}$, $\totalpdv := [\totalpdc{1},\dots,\totalpdc{\pdim}]^{\top}$ and $\vorticitydv := [\vorticitydc{1},\dots,\vorticitydc{\wdim}]^{\top}$. Using matrix notation, \eqref{eq::ns_meevc_weak_form_discrete_expansion} can be expressed more compactly as
		\begin{equation}
			\begin{dcases}
				\text{Find } \velocitydv\in \mathbb{R}^{\udim}, \totalpdv\in \mathbb{R}^{\pdim} \text{ and } \vorticitydv \in \mathbb{R}^{\wdim} \text{ such that:}\\
				\matrixoperator{M} \frac{\mathrm{d}\velocitydv}{\mathrm{d}t} + \matrixoperator{R}\,\velocitydv - \matrixoperator{P}\,\totalpdv = -\nu\,\boldsymbol{l}, \\
				\matrixoperator{N}\frac{\mathrm{d}\vorticitydv}{\mathrm{d}t}  - \frac{1}{2}\matrixoperator{W}\,\vorticitydv + \frac{1}{2}\matrixoperator{W}^{\top}\vorticitydv = \nu\,\matrixoperator{L}\,\vorticitydv, \\
				\matrixoperator{D}\,\velocitydv = 0,
			\end{dcases} \label{eq::ns_meevc_weak_form_discrete_matrix_notation}
		\end{equation}
		The coefficients of the matrices $\matrixoperator{M}$, $\matrixoperator{R}$ and $\matrixoperator{P}$, and the column vector $\boldsymbol{l}$ are given by
		\begin{equation}
			\matrixoperatorc{M}_{ij} := \langle\ubasis{j},\ubasis{i}\rangle_{\Omega}, \quad \matrixoperatorc{R}_{ij} := \langle\vorticityd\times\ubasis{j},\ubasis{i}\rangle_{\Omega}, \quad \matrixoperatorc{P}_{ij} := \langle \pbasis{j},\nabla\cdot\ubasis{i}\rangle_{\Omega}\quad\mathrm{and}\quad l_{i} := \langle\nabla\times\vorticityd,\ubasis{i}\rangle_{\Omega}. \label{eq:matrix_coefficients_1}
		\end{equation}
		Similarly, the coefficients of the matrices $\matrixoperator{N}$, $\matrixoperator{W}$, $\matrixoperator{L}$ and $\matrixoperator{D}$ are given by
		\begin{equation}
			\matrixoperatorc{N}_{ij} := \langle\wbasis{j},\wbasis{i}\rangle_{\Omega}, \quad \matrixoperatorc{W}_{ij} := \langle\wbasis{j},\nabla\cdot\left(\velocityd\,\wbasis{i}\right)\rangle_{\Omega}, \quad \matrixoperator{L}_{ij} := \langle\nabla\times\wbasis{j},\nabla\times\wbasis{i}\rangle_{\Omega} \quad \mathrm{and}\quad  \matrixoperator{D}_{ij} :=\langle\nabla\cdot\ubasis{j},\pbasis{i}\rangle_{\Omega}. \label{eq:matrix_coefficients_2}
		\end{equation}
			
			\subsubsection{Temporal discretization (periodic boundary conditions)}
				Once the spatial discretization is introduced we end up with a set of ordinary differential equations, \eqref{eq::ns_meevc_weak_form_discrete_matrix_notation}. The main objective that drives the choice of the time discretization for \eqref{eq::ns_meevc_weak_form_discrete_matrix_notation} is preservation of invariants: mass, energy, enstrophy, and total vorticity. Not all time integrators satisfy these invariance properties and therefore will spoil all properties obtained so far with the spatial discretization. For the MEEVC scheme the lowest order Gauss time integrator, $s=1$, also known as the \emph{midpoint rule}, is employed. The reason behind this choice has to do with the fact that it enables the construction of a quasi-linear staggered integrator in time. For more details on Gauss time integrators see \cite{Hairer2006}.

\begin{figure}[!ht]
		\centering
		\includegraphics{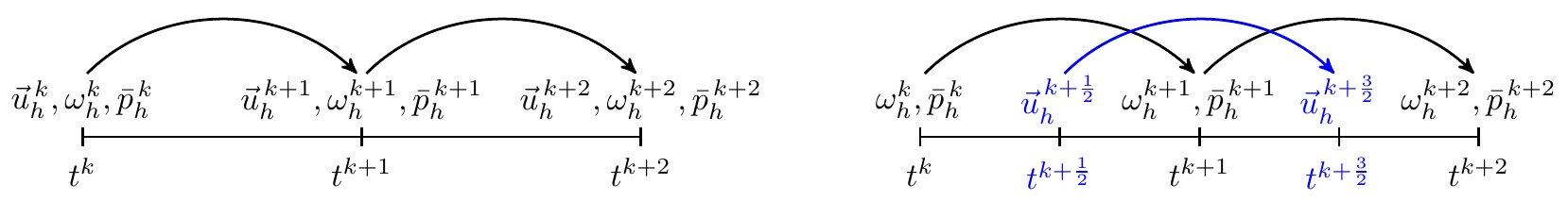}
		\caption{Diagram of the time stepping. Left: all unknowns at the same time instant, as in \eqref{eq::ns_meevc_weak_form_discrete_matrix_notation} Right: staggered in time unknowns, as in \eqref{eq::ns_meevc_weak_form_discrete_staggered_gauss_matrix_notation}.}
		\label{fig:time_stepping}
	\end{figure}

	When applied to the solution of a 1D ordinary differential equation of the form
	\begin{equation}
		\begin{dcases}
			\frac{\mathrm{d}f}{\mathrm{d}t} = g(f(t),t), \\
			f(0) = f_{0},
		\end{dcases}
	\end{equation}
	the one stage Gauss integrator results in the following implicit time stepping scheme
	\begin{equation}
		\frac{f^{k} - f^{k-1}}{\Delta t} = g\left(\frac{f^{k}+f^{k-1}}{2},t+\frac{\Delta t}{2}\right), \quad k=1,\dots,M, \label{eq:gauss_integrator_1D}
	\end{equation}
	where $f^{0} = f_{0}$, $\Delta t$ is the time step and $M$ is the number of time steps. The direct application of \eqref{eq:gauss_integrator_1D} to the discrete weak form \eqref{eq::ns_meevc_weak_form_discrete_matrix_notation} results in a fully nonlinear implicit scheme.This means that the resulting system of equations is a fully coupled set of nonlinear equations, which requires a computationally expensive iterative procedure to solve. To circumvent this penalty, instead of defining all the unknown physical quantities $\velocityd$, $\vorticityd$ and  $\totalpd$, at the same time instants $t^{k}$ we choose to stagger them in time. In this way it is possible to obtain two systems of quasi-linear equations. The unknown vorticity and total pressure are defined at the integer time instants $\vorticitydprevious$, $\totalpdprevious$ and the unknown velocity is defined at the intermediate time instants $\velocitydmidb$, see \figref{fig:time_stepping}. Taking into account this staggered approach, the fully discrete counterpart of \eqref{eq::ns_meevc_weak_form_discrete_matrix_notation} can be rewritten as
	\begin{equation}
			\begin{dcases}
				\text{Find } \velocitydvmidnext\in \mathbb{R}^{\udim}, \totalpdvnext\in \mathbb{R}^{\pdim} \text{ and } \vorticitydvnext \in \mathbb{R}^{\wdim} \text{ such that:}\\
				\matrixoperator{M} \frac{\velocitydvmidnext - \velocitydvmidb}{\Delta t} + \matrixoperator{R}^{k+1}\,\frac{\velocitydvmidnext + \velocitydvmidb}{2} - \matrixoperator{P}\,\totalpdvnext = -\nu\,\boldsymbol{l}^{k+1}, \\
				\matrixoperator{N}\frac{\vorticitydvnext - \vorticitydvprevious}{\Delta t}  - \frac{1}{2}\matrixoperator{W}^{k+\frac{1}{2}}\,\frac{\vorticitydvnext + \vorticitydvprevious}{2} + \frac{1}{2}\left(\matrixoperator{W}^{k+\frac{1}{2}}\right)^{\top}\frac{\vorticitydvnext + \vorticitydvprevious}{2} = \nu\,\matrixoperator{L}\,\frac{\vorticitydvnext + \vorticitydvprevious}{2}, \\
				\matrixoperator{D}\,\velocitydvmidnext = 0,
			\end{dcases} \label{eq::ns_meevc_weak_form_discrete_staggered_gauss_matrix_notation}
	\end{equation}
	where, for compactness of notation, we have set
	\begin{equation}
		\velocitydmida := \frac{\velocitydnext + \velocitydprevious}{2} \quad \mathrm{and}\quad \vorticitydmida := \frac{\vorticitydnext + \vorticitydprevious}{2}. \label{eq:mid_steps_compact_notation}
	\end{equation}
	Note that $\velocitydmidb$ and $\vorticitydprevious$ are known at the start of each time step. All matrix operators are as in \eqref{eq:matrix_coefficients_1} and \eqref{eq:matrix_coefficients_2}, with the exception of $\matrixoperator{R}^{k+1}$, $\matrixoperator{W}^{k+\frac{1}{2}}$ and $\boldsymbol{l}^{k+1}$, the coefficients of which are
	\begin{equation}
		\matrixoperatorc{R}^{k+1}_{ij} := \langle\vorticitydnext\times\ubasis{j},\ubasis{i}\rangle_{\Omega}, \quad l^{k+1}_{i} := \langle\nabla\times\vorticitydnext,\ubasis{i}\rangle_{\Omega}\quad\mathrm{and}\quad \matrixoperatorc{W}^{k+\frac{1}{2}}_{ij} := \langle\wbasis{j},\nabla\cdot\left(\velocitydmidb\,\wbasis{i}\right)\rangle_{\Omega}. \label{eq:matrix_coefficients_staggered}
	\end{equation}
	
	To start the iteration procedure $\velocitydmidstart$ and $\vorticitydstart$ are required. Since only $\velocitydstart$ and $\vorticitydstart$ are known, the first time step needs to be implicit. The remaining time steps can then be computed explicitly with \eqref{eq::ns_meevc_weak_form_discrete_staggered_gauss_matrix_notation}.
	
		\subsubsection{Full discretization (non-periodic boundary conditions)}
		
		% Why can we not impose tangential b.c. strongly on Hdiv-conforming elements?
		
		In general, functions in $\hdiv{\Omega}$-conforming finite element spaces like the Raviart-Thomas elements satisfy weaker continuity requirements than functions contained in $(H^1(\Omega))^2$-conforming spaces, like continuous Galerkin vector elements. In particular, the degrees of freedom of the Raviart-Thomas finite element are calculated in terms of moments of the normal components of the vector along the faces of the cell \cite{kirby2012,brezzi1991}. As a result, one has no control over the tangential component along a face of functions in this space. It follows that tangential boundary conditions cannot be imposed strongly.
		
		This problem was considered in \cite{DeDiego2019a} and three solutions were suggested which consisted in imposing vorticity boundary conditions on the vorticity transport equation. Of these 3 methods, the kinematic Neumann vorticity boundary condition was shown to yield the best approximations and therefore only this method will be considered here. Basically, an additional vorticity $\tilde{\omega}$ is introduced together with the equation $\tilde{\omega} = \mathrm{curl}_h \vec{u}$, where $\mathrm{curl}_h$ is interpreted in the following sense:
		\begin{equation}
		\label{eq:kinematic_relation_weak}
		\tilde{\omega} = \mathrm{curl}_h \vec{u}\iff\langle \tilde{\omega}, \tilde{\xi} \rangle_{\Omega}= \langle \vec{u}, \curl \tilde{\xi} \rangle_{\Omega}- \int_{\partial \Omega} \tilde{\xi} \left(\vec{u}\cdot \vec{\tau} \right)\ \mathrm{d} \Gamma \quad \forall \tilde{\xi} \in \hcurl{\Omega}.
		\end{equation}
		Here, $\vec{\tau}$ is the tangential unit vector along $\partial\Omega$. When solving the vorticity transport equation, the boundary condition $\curl\omega\times\vec{n} = \curl\tilde{\omega}\times\vec{n}$ on $\partial\Omega$ is imposed weakly. Note that although $\omega = \tilde{\omega}$ might hold at the continuous level, at the discrete level this equality does not hold in general. Below, the MEEVC discretization for turbidity currents with kinematic Neumann boundary conditions is described.

		\subsection{MEEVC discretization for turbidity currents}

		Here, the modified MEEVC scheme with kinematic Neumann boundary conditions is adapted for the computation of the lock exchange flow presented in Section \ref{S2:lock_exchange}. The first step is to write \eqref{eq:NS_particle} in the $(\vec{u}, \omega, p)$ formulation of the MEEVC scheme with the additional variable $\tilde{\omega}$:
		\begin{equation}\label{eq:NS_particle_uw}
			\begin{dcases}
				\pddt{\vec{u}} + \omega \times \vec{u} + \gradient \bar{p} = - \frac{1}{\sqrt{\mathrm{Gr}}} \curl \omega + \phi \vec{e}_g,\\
				\frac{\partial \omega}{\partial t} + \frac{1}{2} \left( \vec{u} \cdot \nabla \right) \omega + \frac{1}{2} \nabla \cdot \left( \vec{u} \omega \right)  = \frac{1}{\sqrt{\mathrm{Gr}}} \Delta \omega + \gradient \phi \times \vec{e}_g, \\
				\pddt{\phi} + \frac{1}{2} \left( \vec{u}_p \cdot \gradient \right) \phi + \frac{1}{2} \divergence \left( \vec{u}_p \phi \right) = \frac{1}{\sqrt{\mathrm{Gr} \mathrm{Sc}^2}} \Delta \phi, \\
				\tilde{\omega} = \mathrm{curl}_h\vec{u},\\
				\divergence \vec{u} = 0.
			\end{dcases}
		\end{equation}
		The system \eqref{eq:NS_particle_uw} is considered with the boundary conditions presented in Section \ref{S2:lock_exchange} (see Figure \ref{fig:bc}) and the additional vorticity constraint (which may be interpreted as a boundary condition for vorticity evolution equation)
		\begin{equation}
			\curl\omega\times\vec{n} = \curl\tilde{\omega}\times\vec{n} \quad \text{on $\Gamma_1\cup\Gamma_3$}.
		\end{equation}
		In \eqref{eq:NS_particle_uw}, $\vec{u}_p = \vec{u} + u_s \vec{e}_g$ and the identity $\curl (\phi \vec{e}_g) = \gradient \phi \times \vec{e}_g$ have been used. The convective term for the particle transport equation has been written in skew-symmetric form to conserve the quadratic mean $\langle \phi, \phi \rangle$ at the discrete level. 

		The next step is the construction of a weak formulation. The velocity $\vec{u}$ is sought in the function space $H_0 (\mathrm{div}, \Omega)$, defined by
		\begin{equation}
		\label{eq:H0_div}
		H_0( \textrm{div}, \Omega) = \left\lbrace \vec{u} \in H(\textrm{div}, \Omega)\ |\ \vec{u} \cdot \vec{n} = 0 \ \text{on $\partial \Omega$} \right\rbrace ,
		\end{equation}		
		such that $\vec{u}\cdot\vec{n} = 0$ is strongly imposed along $\partial \Omega$. On the other hand, the vorticity $\omega$ is assumed to be contained in the space $H_{(2,4)} (\mathrm{curl}, \Omega)$, defined
		\begin{equation}\label{eq:hcurl_024}
			H_{(2,4)} (\mathrm{curl}, \Omega) = \left\lbrace \omega \in H(\textrm{curl}, \Omega)\ |\ \omega = 0\ \text{on $\Gamma_2 \cup \Gamma_4$} \right\rbrace .
		\end{equation}

		The particle concentration $\phi$ is contained in the space of functions with square integrable gradients, $H(\mathrm{grad},\Omega)$. The boundary conditions, given by \eqref{eq:p_no_flux_bc2} and \eqref{eq:p_flux_bc}, must be imposed weakly by means of Neumann boundary conditions. After integrating by parts, the weak form of the particle transport equation is,
		\begin{equation}\label{eq:prt_transport_weak}
			\begin{split}
				\langle \pddt{\phi}, \zeta \rangle_{\Omega}+ \frac{1}{2} \langle \divergence \left( \vec{u}_p \phi \right), \zeta \rangle_{\Omega}- \frac{1}{2}\langle \phi, \divergence \left( \vec{u}_p \zeta \right) \rangle_{\Omega}+ \frac{1}{2} \int_{\partial \Omega} \zeta \phi\left( \vec{u}_p \cdot \vec{n}\right)\ \mathrm{d}\Gamma = \\ 
				= \frac{1}{\sqrt{\mathrm{Gr} \mathrm{Sc}^2}} \left[ - \langle \gradient \phi, \gradient \zeta \rangle_{\Omega}+ \int_{\partial \Omega} \zeta \left( \gradient \phi \cdot \vec{n} \right)\ \mathrm{d}\Gamma \right], \quad \forall \zeta \in \hgrad{\Omega}. \\
			\end{split}
		\end{equation}

		The term $\gradient \phi \cdot \vec{n}$ is equal to zero along $\Gamma_2 \cup \Gamma_3 \cup \Gamma_4$. On the other hand, $\vec{u}_p \cdot \vec{n} = u_s \vec{e}_g \cdot \vec{n}$ will be zero along the lateral faces, $\Gamma_2 \cup \Gamma_4$. Finally, $\gradient \phi \cdot \vec{n} = - \sqrt{\mathrm{Gr} \mathrm{Sc}^2} \phi u_s (\vec{e}_g \cdot \vec{n})$ on $\Gamma_1$. These relations imply that \eqref{eq:prt_transport_weak} can be written as follows
		\begin{equation}\label{eq:prt_transport_weak2}
			\begin{split}
				& \langle \pddt{\phi}, \zeta \rangle_{\Omega}+ \frac{1}{2} \langle \divergence \left( \vec{u}_p \phi \right), \zeta \rangle_{\Omega}- \frac{1}{2}\langle \phi, \divergence \left( \vec{u}_p \zeta \right) \rangle_{\Omega}= \\ 
				= - & \frac{1}{\sqrt{\mathrm{Gr} \mathrm{Sc}^2}} \langle \gradient \phi, \gradient \zeta \rangle_{\Omega}+ \int_{\Gamma_1} \zeta \phi u_s \left( \vec{e}_g \cdot \vec{n} \right)\ \mathrm{d}\Gamma - \frac{1}{2} \int_{\Gamma_1 \cup \Gamma_3} \zeta \phi u_s \left( \vec{e}_g \cdot \vec{n}\right)\ \mathrm{d}\Gamma = \\
				= - & \frac{1}{\sqrt{\mathrm{Gr} \mathrm{Sc}^2}} \langle \gradient \phi, \gradient \zeta \rangle_{\Omega}+ \frac{1}{2} \int_{\Gamma_1} \zeta \phi u_s\ \mathrm{d}\Gamma - \frac{1}{2} \int_{\Gamma_3} \zeta \phi u_s\ \mathrm{d}\Gamma, \quad \forall \zeta \in \hgrad{\Omega}.
			\end{split}
		\end{equation}

		The weak form of \eqref{eq:NS_particle_uw} can then be written as
		\begin{equation}
			\label{eq:NS_prt_weak_1}
				\begin{dcases}
					\text{Find $\vec{u} \in H_0( \textrm{div}, \Omega),\ p \in L^2(\Omega),\ \omega \in H_{(2,4)} (\mathrm{curl}, \Omega)$, $\phi \in \hgrad{\Omega}$} &\\
					\text{and $\tilde{\omega} \in H_{(2,4)} (\mathrm{curl})$ such that} & \\
					\langle \frac{\partial\vec{u}}{\partial t}, \vec{v} \rangle_{\Omega}+ \langle \omega \times \vec{u} , \vec{v} \rangle_{\Omega}- \langle \bar{p}, \nabla \cdot \vec{v} \rangle_{\Omega}= -\frac{1}{\sqrt{\mathrm{Gr}}} \langle \nabla \times \omega, \vec{v} \rangle_{\Omega}+ \langle \phi \vec{e}_g, \vec{v} \rangle_\Omega, & \forall \vec{v} \in H_0( \textrm{div}, \Omega), \\
					\langle \frac{\partial \omega}{\partial t}, \xi \rangle_{\Omega}+C_{\vec{u}}(\omega,\xi) = \frac{1}{\sqrt{\mathrm{Gr}}} \left[ -\langle \nabla \times \omega, \nabla \times \xi \rangle_{\Omega}  + \langle \xi,  \left( \curl \tilde{\omega} \right) \times \vec{n} \rangle_{\partial\Omega}\ \right]+ \langle\gradient \phi \times \vec{e}_g, \xi \rangle_{\Omega}, &\forall \xi \in H_{(2,4)} (\mathrm{curl}, \Omega), \\
					\langle \pddt{\phi}, \zeta \rangle_{\Omega}+ D_{\vec{u}}(\phi,\zeta) = - \frac{1}{\sqrt{\mathrm{Gr} \mathrm{Sc}^2}} \langle \gradient \phi, \gradient \zeta \rangle_\Omega, & \forall \zeta \in \hgrad{\Omega} , \\
					\langle \tilde{\omega}, \tilde{\xi} \rangle_{\Omega}= \langle \vec{u}, \curl \tilde{\xi} \rangle_{\Omega}- \langle \tilde{\xi},\vec{u}\cdot \vec{\tau} \rangle_{\partial \Omega} & \forall \tilde{\xi} \in \hcurl{\Omega} \\
					\langle \nabla \cdot \vec{u}, q \rangle_{\Omega}= 0, & \forall q \in L^2(\Omega).
				\end{dcases}
			\end{equation}
			where we have introduced the velocity dependent bilinear forms $C_{\vec{u}}: H_{(2,4)}(\mathrm{curl}, \Omega)\times H_{(2,4)}(\mathrm{curl}, \Omega) \to \mathbb{R}$ and $D_{\vec{u}}: H(\mathrm{grad}, \Omega)\times H(\mathrm{grad}, \Omega) \to \mathbb{R}$ defined by
			\begin{equation}
				C_{\vec{u}}(\omega,\xi) = \frac{1}{2} \langle \divergence \left( \vec{u} \omega \right), \xi \rangle_{\Omega}- \frac{1}{2}\langle \omega, \divergence \left( \vec{u} \xi \right) \rangle_\Omega
			\end{equation}
			and
			\begin{equation}
				D_{\vec{u}}(\phi,\zeta;\vec{u}) = \frac{1}{2} \langle \divergence \left( \vec{u}_p \phi \right), \zeta \rangle_{\Omega}- \frac{1}{2}\langle \phi, \divergence \left( \vec{u}_p \zeta \right) \rangle_{\Omega}- \frac{1}{2} \langle \zeta, \phi u_s\rangle_{\Gamma_{1}} + \frac{1}{2} \langle \zeta, \phi u_s\rangle_{\Gamma_{3}}\, ,
			\end{equation}
			where the convective particle velocity is $\vec{u}_p = \vec{u} + u_s \vec{e}_g$.

%			\begin{equation}
%				D_{\omega}(\omega, \xi) = - \frac{1}{\sqrt{\mathrm{Gr}}} \langle \nabla \times \omega, \nabla \times \xi \rangle_{\Omega}  + \frac{1}{\sqrt{\mathrm{Gr}}}\langle \xi,  \left( \curl \tilde{\omega} \right) \times \vec{n} \rangle_{\partial\Omega}\,,
%			\end{equation}
%			and $S_{\phi}: \phi \otimes \zeta \in H(\mathrm{grad}, \Omega) \otimes H(\mathrm{grad}, \Omega) \mapsto S_{\phi}(\phi, \zeta)$ with
%			\begin{equation}
%				S_{\phi}(\phi, \zeta) =  \frac{1}{2} \langle \zeta, \phi u_s\rangle_{\Gamma_{1}} - \frac{1}{2} \langle \zeta, \phi u_s\rangle_{\Gamma_{2}}\,.
%			\end{equation}

			The finite element discretization of \eqref{eq:NS_prt_weak_1} can be developed once a set of finite element spaces is chosen for the variables, these are extended in terms of the basis functions and a time integrator is implemented. Following the steps taken in the MEEVC scheme with kinematic Neumann boundary conditions and taking into account that the space $\mathrm{CG}_N$ is chosen for $\phi_h$, the following algorithm is developed:

			\noindent Given $\vec{u}_h^{\,k+\frac{1}{2}}$, $\omega_h^k$ and $\phi^k$,\\

			\noindent \textbf{Step 1}. Find $\tilde{\omega}_h^{k+1}\in \mathrm{CG}_{N,(2,4)} = \left\lbrace \omega_h \in \mathrm{CG}_{N}\ |\ \omega_h = 0\ \text{on $\Gamma_2 \cup \Gamma_4$} \right\rbrace$ such that:
			\begin{equation}\label{eq:step_tilde_omega}
				\langle \tilde{\omega}_h^{k+\frac{1}{2}}, \xi_h \rangle_{\Omega}= \langle \vec{u}_h^{\,k+\frac{1}{2}}, \curl \xi_h \rangle_\Omega, \quad \forall \xi_h \in \mathrm{CG}_{N,(2,4)}.
			\end{equation}

			\noindent \textbf{Step 2}. Find $\phi_h^{k+1}\in \mathrm{CG}_N$ such that:
			\begin{equation}\label{eq:step_phi}
					\langle \frac{\phi^{k+1}_h - \phi^{k}_h}{\Delta t}, \zeta_h \rangle_{\Omega}+D_{\vec{u}_h^{k+\frac{1}{2}}}\left(\frac{\phi^{k+1}_h + \phi^{k}_h}{2},\zeta_h \right)= - \frac{1}{\sqrt{\mathrm{Gr} \mathrm{Sc}^2}} \langle \gradient \left( \frac{\phi^{k+1}_h + \phi^{k}_h}{2} \right) , \gradient \zeta_h \rangle_{\Omega}\quad  \forall \zeta_h \in \mathrm{CG}_N
			\end{equation}

			\noindent \textbf{Step 3}. Find $\omega_h^{k+1}\in \mathrm{CG}_{N,(2,4)}$ such that:
				\begin{equation}\label{eq:step_omega}
					\begin{split}
						\langle \frac{\omega^{k+1}_h - \omega^{k}_h}{\Delta t}, \xi_h \rangle_{\Omega} + C_{\vec{u}^{k+\frac{1}{2}}_h}\left(\frac{\omega^{k+1}_h + \omega^{k}_h}{2},\xi_h \right)= \langle\gradient \tilde{\phi}_h^{k+\frac{1}{2}} \times \vec{e}_g, \xi_h \rangle_\Omega + \\
						+ \frac{1}{\sqrt{\mathrm{Gr}}} \left[ - \langle \nabla \times \omega^{k+1}_h, \vec{v}_h \rangle_{\Omega}+ \langle \xi_h, \left( \curl \tilde{\omega}_h^{k+\frac{1}{2}} \right) \times \vec{n} \rangle_{\partial \Omega} \right],\quad \forall \xi_h \in \mathrm{CG}_{N,(2,4)},
					\end{split}
				\end{equation}
				\noindent with
				\begin{equation}\label{eq:phi_tilde}
					\tilde{\phi}_h^{k+\frac{1}{2}} = \frac{\phi^{k+1}_h + \phi^{k}_h}{2}.
				\end{equation} 

				\noindent \textbf{Step 4}. Find $\left( \vec{u}^{\,k+\frac{3}{2}}_h,\ \bar{p}^{k+1} \right) \in \left( \mathrm{RT}_{N,0},\ \mathrm{DG}_{N-1} \right) $ such that:
				\begin{equation}\label{eq:step_mixed}
					\begin{dcases}
						\langle \frac{\vec{u}^{\,k+\frac{3}{2}}_h - \vec{u}^{\,k+\frac{1}{2}}_h}{\Delta t}, \vec{v}_h \rangle_{\Omega}+ \langle \omega^{k+1}_h \times \frac{\vec{u}^{\,k+\frac{3}{2}}_h + \vec{u}^{\,k+\frac{1}{2}}_h}{2} , \vec{v}_h \rangle_{\Omega}- \langle \bar{p}^{k+1}_h, \nabla \cdot \vec{v}_h \rangle_\Omega=\\
						 -\frac{1}{\sqrt{\mathrm{Gr}}} \langle \nabla \times \omega^{k+1}_h, \vec{v}_h \rangle_{\Omega}+ \langle \phi_h^{k+1} \vec{e}_g, \vec{v} \rangle_{\Omega},& \forall \vec{v}_h \in\mathrm{RT}_{(N,0)}, \\[0.6em]
						\langle \nabla \cdot \vec{u}^{\,k+\frac{3}{2}}_h, q_h \rangle_{\Omega}= 0,& \forall q_h \in \mathrm{DG}_{N-1}.
					\end{dcases}
				\end{equation}

	\section{Energy balance}\label{S2:eb}
		The progressive conversion of potential energy into kinetic energy is the main mechanism that drives a turbidity current. The fluid motion decays over time due to viscous dissipation and the gradual loss of mass due to sedimentation. In this section, the energy balance equations are derived at the continuous level in the way of \cite{necker2002,necker2005,parkinson2014,espath2014}. Subsequently, the degree to which the discrete counterpart of the energy balance holds is analyzed.

		\subsection{Conservation of energy at the continuous level}\label{ss2:cons_E_sec2}
			Denoting the potential energy by $E_p$, the kinetic energy by $\mathcal{K}$ and the viscous and sedimentation dissipation rates with $\varepsilon_v$ and $\varepsilon_s$, the global energy equation of the flow is given by,
			\begin{equation}\label{eq:energy_particle}
				\ddt{} \left( \mathcal{K} + E_p \right) = - \varepsilon_v - \varepsilon_s.
			\end{equation}

			An expression for the conservation of kinetic energy can be derived from the momentum equation in its weak formulation by specifying $\vec{v} = \vec{u}$, see \eqref{eq:NS_prt_weak_1}, such that,
			\begin{equation}\label{eq:energy_eq}
				\ddt{\mathcal{K}} = -\frac{1}{\sqrt{\mathrm{Gr}}} \langle \nabla \times \omega, \vec{u} \rangle_{\Omega}+ \langle \phi \vec{e}_g, \vec{u} \rangle_\Omega.
			\end{equation}

			The potential energy, in its non-dimensional form, is defined,
			\begin{equation}\label{eq:E_p}
				E_p = \int_\Omega \phi y\ \mathrm{d}\Omega,
			\end{equation} 
		    such that its time derivative is given by
			\begin{equation}
				\begin{split}
				\ddt{E_p} = \int_\Omega \frac{\partial\phi}{\partial t} y\ \mathrm{d}\Omega \stackrel{\eqref{eq:NS_prt_weak_1}}{=} -\frac{1}{2} \langle \divergence \left( \vec{u}_p \phi \right), y \rangle_{\Omega} + \frac{1}{2}\langle \phi, \divergence \left( \vec{u}_p y \right) \rangle_{\Omega} + \\
				+ \frac{1}{\sqrt{\mathrm{Gr} \mathrm{Sc}^2}} \langle \gradient \phi, \vec{e}_{g} \rangle_{\Omega}+\frac{1}{2} \langle \zeta, \phi u_s\rangle_{\Gamma_{1}} - \frac{1}{2} \langle \zeta, \phi u_s\rangle_{\Gamma_{3}}.
				\end{split}
			\end{equation} 
			Expanding the terms of this expression and integrating by parts yields 
			\begin{equation}\label{eq:dt_E_p}
				\ddt{E_p} = \int_\Omega \frac{\partial\phi}{\partial t} y\ \mathrm{d}\Omega = -\langle \phi , \vec{u}_{p}\cdot\vec{e}_{g} \rangle_{\Omega} + \frac{1}{\sqrt{\mathrm{Gr} \mathrm{Sc}^2}} \langle \gradient \phi, \vec{e}_{g} \rangle_{\Omega}.
			\end{equation} 

			If we now combine \eqref{eq:dt_E_p} with \eqref{eq:energy_eq}  we obtain the following expression for the energy budget equation 
			\begin{equation}\label{eq:energy_eq2}
				\ddt{\mathcal{K}} + \ddt{E_p} =  -\frac{1}{\sqrt{\mathrm{Gr}}} \langle \nabla \times \omega, \vec{u} \rangle_{\Omega} - u_{s}\,\langle \phi , 1 \rangle_{\Omega} + \frac{1}{\sqrt{\mathrm{Gr} \mathrm{Sc}^2}} \langle \gradient \phi, \vec{e}_{g} \rangle_{\Omega}.
			\end{equation}

			This equation yields an expression for the dissipation terms due to viscosity, denoted $\epsilon_v$, and sedimentation, $\epsilon_s$, such that,
			\begin{align}
				\label{eq:eps_v} & \epsilon_v = \frac{1}{\sqrt{\mathrm{Gr}}} \langle \nabla \times \omega, \vec{u} \rangle_{\Omega}, \\
				\label{eq:eps_s0} & \epsilon_s = u_{s}\,\langle \phi , 1 \rangle_{\Omega} - \frac{1}{\sqrt{\mathrm{Gr} \mathrm{Sc}^2}} \langle \gradient \phi, \vec{e}_{g} \rangle_{\Omega}.
			\end{align}

			The energy balance of a turbidity current is generally written in terms of time integrated values \cite{necker2002,necker2005,espath2014,parkinson2014}. Equation \eqref{eq:energy_eq} can be written as follows
			\begin{equation}\label{eq:energy_eq_3}
				\mathcal{K} + E_p + E_v + E_s = \mathcal{K}_0 + E_{p,0},
			\end{equation}

			\noindent where $\mathcal{K}_0$ and $E_{p,0}$ are the initial kinetic and potential energy and $E_v$ and $E_s$ are the total energy losses to viscous dissipation and suspended particles. $E_v$ and $E_s$ are given by
			\begin{equation}\label{eq:eps_s_v}
				E_v = \int_0^t \epsilon_v\ \mathrm{d}\tau \quad \text{and} \quad E_s = \int_0^t \epsilon_s\ \mathrm{d}\tau.
			\end{equation}

%ANALYZE ENERGY CONSERVATION PROPERTIES OF THE ADAPTED MEEVC SCHEME
			\subsection{Conservation of energy at the discrete level}
				A discrete energy balance equation, equivalent to \eqref{eq:energy_eq_3}, can be derived from the discretized equations of motion, given by \eqref{eq:step_tilde_omega} to \eqref{eq:step_mixed}. The discretized energy balance equation is constructed by imitating the derivation carried out in Section \ref{ss2:cons_E_sec2} and taking into account which steps and vector identities are exact at the discrete level.

				A discrete statement of the conservation of kinetic energy can be written by setting $\vec{v}_h = \frac{1}{2} \vec{u}^{\,k+\frac{3}{2}}_h+ \frac{1}{2} \vec{u}^{\,k+\frac{1}{2}}_h$ in the discrete momentum equation. The following equation can be derived,
				\begin{equation}\label{eq:disc_energy_eq}
					\frac{1}{2} \langle \vec{u}^{\,k+\frac{3}{2}}_h , \vec{u}^{\,k+\frac{3}{2}}_h \rangle_{\Omega}- \frac{1}{2} \langle \vec{u}^{\,k+\frac{1}{2}}_h , \vec{u}^{\,k+\frac{1}{2}}_h \rangle_\Omega= - \frac{\Delta t}{\sqrt{\mathrm{Gr}}} \langle \nabla \times \omega_h^{k+1}, \frac{\vec{u}^{\,k+\frac{3}{2}}_h + \vec{u}^{\,k+\frac{1}{2}}_h}{2} \rangle_{\Omega}+ \Delta t \langle \phi^{k+1} \vec{e}_g, \frac{\vec{u}^{\,k+\frac{3}{2}}_h + \vec{u}^{\,k+\frac{1}{2}}_h}{2} \rangle_\Omega.
				\end{equation}

				The discrete kinetic energy at the time step $t^k$ is given by $\mathcal{K}_h^k = 1/2 \langle \velocitydprevious, \velocitydprevious \rangle$. The discrete counterpart of $\epsilon_v$ is defined as
				\begin{equation}\label{eq:eps_v_discrete}
					\epsilon^{k+1}_{v,h} = \frac{1}{\sqrt{\mathrm{Gr}}} \langle \nabla \times \omega_h^{k+1}, \frac{\vec{u}^{\,k+\frac{3}{2}}_h + \vec{u}^{\,k+\frac{1}{2}}_h}{2} \rangle_\Omega.
				\end{equation}

				Therefore, Equation \eqref{eq:disc_energy_eq} can be written as follows,
				\begin{equation}\label{eq:disc_energy_eq2}
					\mathcal{K}^{k+1}_h - \mathcal{K}^{k}_h= - \Delta t\ \epsilon^{k+1}_{v,h} + \Delta t \langle \phi^{k+1} \vec{e}_g, \frac{\vec{u}^{\,k+\frac{3}{2}}_h + \vec{u}^{\,k+\frac{1}{2}}_h}{2} \rangle_\Omega.
				\end{equation}

				The potential energy and the dissipation due to sediementation can be introduced by using the discrete transport equation for the particle phase, Equation \eqref{eq:step_phi}. To this end, the discrete test function is set to $\zeta_h = y$ and the convective terms are rewritten as follows,
				\begin{equation}\label{eq:disc_conv}
					\frac{1}{2} \langle \divergence \left( \vec{u}_h^{\,k+\frac{1}{2}} \frac{\phi^{k+1}_h + \phi^{k}_h}{2} \right), y \rangle_{\Omega}- \frac{1}{2}\langle \frac{\phi^{k+1}_h + \phi^{k}_h}{2}, \divergence \left( \vec{u}_h^{\,k+\frac{1}{2}} y \right) \rangle_{\Omega}= \langle \tilde{\phi}^{k+\frac{1}{2}}, \vec{u}_h^{\,k+\frac{1}{2}} \cdot \vec{e}_g \rangle_\Omega,
				\end{equation}
				\begin{equation}\label{eq:disc_conv2}
					\frac{1}{2} \langle \divergence \left( u_s\vec{e}_g \frac{\phi^{k+1}_h + \phi^{k}_h}{2} \right), y \rangle_{\Omega}- \frac{1}{2}\langle \frac{\phi^{k+1}_h + \phi^{k}_h}{2}, \divergence \left( u_s\vec{e}_g y \right) \rangle_{\Omega}+ \frac{1}{2} \int_{\Gamma_1} y \frac{\phi^{k+1}_h + \phi^{k}_h}{2} u_s (\vec{e}_g \cdot \vec{n} )\ \mathrm{d}\Gamma = - u_s \langle \tilde{\phi}_h^{k+\frac{1}{2}}, 1 \rangle_\Omega,
				\end{equation}

				\noindent where integration by parts has been applied on the first terms in both cases, an operation which is exact at the discrete level because $\divergence \vec{u}_h = 0$ and $\vec{u}_h \cdot \vec{n} = 0$. For simplicity of notation, $\tilde{\phi}^{k+\frac{1}{2}}$ is used, see \eqref{eq:phi_tilde}. The discrete counterparts of $E_p$ and $\epsilon_s$ are defined as,
				\begin{equation}\label{eq:Ep_disc}
					E_{p,h}^k = \langle \phi_h^k, y \rangle_\Omega,
				\end{equation}
				\begin{equation} \label{eq:eps_s_disc}
					\epsilon_{s,h}^k =  u_s \langle \tilde{\phi}_h^k, 1 \rangle_{\Omega}- \frac{1}{\sqrt{\mathrm{Gr}\mathrm{Sc}^2}}\langle \gradient \tilde{\phi}_h^k, \vec{e}_g \rangle_\Omega,
				\end{equation}

				\noindent such that \eqref{eq:step_phi} can be rewritten as follows,
				\begin{equation}\label{eq:transport_discrete}
					E_{p,h}^{k+1} -  E_{p,h}^k = - \Delta t\ \epsilon_{s,h}^{k+\frac{1}{2}} - \Delta t \langle \frac{\phi_h^{k+1} + \phi_h^k}{2}\vec{e}_g, \vec{u}_h^{\,k+\frac{1}{2}} \rangle_\Omega.
				\end{equation}

				A discrete energy balance equation analogous to \eqref{eq:energy_particle}  is obtained by summing \eqref{eq:disc_energy_eq2} and \eqref{eq:transport_discrete},
				\begin{equation}\label{eq:energy_balance_2}
					\mathcal{K}_h^{k+\frac{3}{2}} + E_{p,h}^{k+1} = \mathcal{K}_h^{k+\frac{1}{2}} + E_{p,h}^k - \Delta t \left( \epsilon_{v,h}^{k+1} + \epsilon_{s,h}^{k+\frac{1}{2}} \right)  + \frac{\Delta t}{2} \left( \langle \phi_h^{k+1}\vec{e}_g, \vec{u}_h^{\,k+\frac{3}{2}} \rangle_{\Omega}- \langle \phi_h^{k}\vec{e}_g, \vec{u}_h^{\,k+\frac{1}{2}} \rangle_{\Omega}\right).
				\end{equation}

				The last term in  \eqref{eq:energy_balance_2} is a residual term due to the staggering in time of the momentum equation with respect to particle transport equation. The temporal evolution of the residual can be understood by deriving the integral form of \eqref{eq:energy_balance_2}. The following expressions are obtained at the time iterations $k,\ k-1,\  ...\ ,1, 0$ :
				\begin{align*}
					& \mathcal{K}_h^{k+\frac{3}{2}} + E_{p,h}^{k+1} = \mathcal{K}_h^{k+\frac{1}{2}} + E_{p,h}^k - \Delta t \left( \epsilon_{v,h}^{k+1} + \epsilon_{s,h}^{k+\frac{1}{2}} \right)  + \frac{\Delta t}{2} \left( \langle \phi_h^{k+1}\vec{e}_g, \vec{u}_h^{\,k+\frac{3}{2}} \rangle_{\Omega}- \langle \phi_h^{k}\vec{e}_g, \vec{u}_h^{\,k+\frac{1}{2}} \rangle_{\Omega}\right), \\
					& \mathcal{K}_h^{k+\frac{1}{2}} + E_{p,h}^{k} = \mathcal{K}_h^{k-\frac{1}{2}} + E_{p,h}^{k-1} - \Delta t \left( \epsilon_{v,h}^{k} + \epsilon_{s,h}^{k-\frac{1}{2}} \right)  + \frac{\Delta t}{2} \left( \langle \phi_h^{k}\vec{e}_g, \vec{u}_h^{\,k+\frac{1}{2}} \rangle_{\Omega}- \langle \phi_h^{k-1}\vec{e}_g, \vec{u}_h^{\,k-\frac{1}{2}} \rangle_{\Omega}\right), \\
					& \quad \vdots \\
					& \mathcal{K}_h^{\frac{5}{2}} + E_{p,h}^{2} = \mathcal{K}_h^{\frac{3}{2}} + E_{p,h}^1 - \Delta t \left( \epsilon_{v,h}^{2} + \epsilon_{s,h}^{\frac{3}{2}} \right)  + \frac{\Delta t}{2} \left( \langle \phi_h^{2}\vec{e}_g, \vec{u}_h^{\,\frac{5}{2}} \rangle_{\Omega}- \langle \phi_h^{1}\vec{e}_g, \vec{u}_h^{\,\frac{3}{2}} \rangle_{\Omega}\right),\\
					& \mathcal{K}_h^{\frac{3}{2}} + E_{p,h}^{1} = \mathcal{K}_h^{\frac{1}{2}} + E_{p,h}^{0} - \Delta t \left( \epsilon_{v,h}^{1} + \epsilon_{s,h}^{\frac{1}{2}} \right)  + \frac{\Delta t}{2} \left( \langle \phi_h^{1}\vec{e}_g, \vec{u}_h^{\,\frac{3}{2}} \rangle_{\Omega}- \langle \phi_h^{0}\vec{e}_g, \vec{u}_h^{\,\frac{1}{2}} \rangle_{\Omega}\right).
				\end{align*}

% & \mathcal{K}_h^1 = \mathcal{K}_h^0 - \Delta t\ \epsilon_{v,h}^1  + \frac{\Delta t}{2} \left( \langle \phi_h^{1}\vec{e}_g, \vec{u}_h^{\,\frac{3}{2}} \rangle + \langle \phi_h^{1}\vec{e}_g, \vec{u}_h^{\,\frac{1}{2}} \rangle \right),

				The summation of all the energy balance equations from $t^{k+1}$ to $t^{1}$ produces the following result,
				\begin{equation}\label{eq:energy_balance_3}
					\mathcal{K}_h^{k+\frac{3}{2}} + E_{p,h}^{k+1} = \mathcal{K}_h^{\frac{1}{2}} + E_{p,h}^{0} - \Delta t \sum_{i=1}^{k+1} \left( \epsilon_{v,h}^i + \epsilon_{v,h}^{i-\frac{1}{2}} \right) + \frac{\Delta t}{2} \left( \langle \phi_h^{k+1}\vec{e}_g, \vec{u}_h^{\,k + \frac{3}{2}} \rangle_{\Omega}- \langle \phi_h^{0}\vec{e}_g, \vec{u}_h^{\,\frac{1}{2}} \rangle_{\Omega} \right).
				\end{equation}

				Given the discrete counterpart of the total dissipated energy,
				\begin{equation}\label{eq:eps_s_v_disc}
					E_{v,h}^{k+1} = \Delta t \sum_{i=1}^{k+1} \epsilon_{v,h}^i \quad \text{and} \quad E^{k+\frac{1}{2}}_{s,h} = \Delta t \sum_{i=1}^{k+1} \epsilon_{s,h}^{i-\frac{1}{2}},
				\end{equation}

				\noindent an integral energy balance equation for the discrete system can be derived,
				\begin{equation}\label{eq:energy_balance_4}
					\mathcal{K}_h^{k+\frac{3}{2}} + E_{p,h}^{k+1} +E_{v,h}^{k+1} + E^{k+\frac{1}{2}}_{s,h} = \mathcal{K}_h^{\frac{1}{2}} + E_{p,h}^{0}  + \frac{\Delta t}{2} \left( \langle \phi_h^{k+1}\vec{e}_g, \vec{u}_h^{\,k + \frac{3}{2}} \rangle_{\Omega}- \langle \phi_h^{0}\vec{e}_g, \vec{u}_h^{\,\frac{1}{2}} \rangle_{\Omega} \right).
				\end{equation}

				Equation \eqref{eq:energy_balance_4} indicates that, due to staggering in time, a discrete statement for the conservation of energy cannot be exact because the different components of the energy are known at different time instants. The residual term of \eqref{eq:energy_balance_4} is not a cumulative error that is unbounded in time. It quantifies the level of mismatch between the staggered variables in time and is proportional only to $\Delta t$.
		
	\section{Numerical results}\label{S2:nr}
		In this section, results are shown for the lock exchange test case described in Section \ref{S2:lock_exchange}. The objective of this section is to validate the modified MEEVC scheme and its robustness. In order to do so, 3 simulations are carried out with different levels of refinement and the results are compared with those given by Parkinson et al. \cite{parkinson2014} and Espath et al. \cite{espath2014}.

		In these computations, the settling velocity is set to $u_{s} = 0.02$, the Grashof and the Schmidt numbers are set to $\mathrm{Gr}=5\times10^6$ ( $\mathrm{Re} \approx 2236$ ) and $\mathrm{Sc}=1$. The simulations presented in the literature compute the flow problem up to $t=60$ for a domain of length $L=19$. The references used for the comparison, denoted Ref. 1 and Ref. 2 hereafter, are briefly explained below.
		\begin{itemize}
			\item \textbf{Ref. 1}, Parkinson et al. \cite{parkinson2014}. A mixed finite element method based on discontinuous Galerkin (DG) elements is used, with linear elements for $\vec{u}$ and $\phi$ and quadratic elements for $p$. A Crank-Nicolson time discretization is used and the resulting non-linear system of equations is solved using two Picard iterations. Fixed and highly adapted meshes are used. Unlike the boundary conditions presented in Section \ref{subsec:le_bc}, $\phi = 0$ is prescribed along the top boundary, resulting in a loss of $<1\%$ of mass and $\sim 3 \%$ of the total energy.
			\item \textbf{Ref. 2}, Espath et al. \cite{espath2014}. A compact sixth-order finite difference scheme is used for the spatial discretization and a third order Adams-Bashforth scheme for time integration. 
		\end{itemize}

		In total, 3 simulations are carried out with the modified MEEVC scheme; these are denoted by Sim 1,2 and 3 and the discretization parameters can be found in Table \ref{tab:le_sim}. The simulations are computed from $t=0$ to $t = 12$ and in a domain of length $L = 13$. At $t=12$, the particle flow is far from reaching the opposite wall and the difference in domain length with the references does not modify the results. H{\"a}rtel et al. \cite{hartel1997} show that the propagating front remains unaffected by the wall up to a distance of $2H$, a condition satisfied in these computations.

		In all of the simulations, the time step is set to $\Delta t = 10^{-3}$, following Espath et al. \cite{espath2014}. An ``equivalent number of cells'' is calculated in order to facilitate comparisons, defined
		\begin{equation}\label{eq:eq_num_cells}
			\text{eq.\ num.\ cells} = \text{num.\ cells} \times \frac{19}{13} \times N^2,
		\end{equation}

		\noindent such that the domain is extended to $L=19$ and the polynomial order is accounted for by triangulating each element into $N^2$ cells. The meshes used in the simulations are depicted in Figure \ref{fig:colorplot_le}.

		\begin{table}[b]
			\centering
				\bgroup
					\def\arraystretch{1.25}
					\caption{Discretization parameters of the three simulations carried out with the modified MEEVC scheme. The Grashof and the Schmidt numbers are set to $\mathrm{Gr}=5\times10^6$ and $\mathrm{Sc}=1$.}
					\label{tab:le_sim}
					\begin{tabular}{cccccccc}
						\hline                                                                                 
						id & $N$ & num. cells & eq. num. cells & $\Delta t$ & $d_W$ & $d_U$ & $d_Q$ \\
						\hline
						Sim 1 & 4 & 1116 & $2.6 \times 10^4$ & $10^{-3}$ & 9169 & 20328 & 11160 \\
						Sim 2 & 4 & 2399 & $5.6 \times 10^4$ & $10^{-3}$ & 19587 & 43576 & 23990 \\
						Sim 3 & 4 & 3599 & $8.4 \times 10^4$ & $10^{-3}$ & 29283 & 65272 & 35990 \\
						\hline
					\end{tabular}
				\egroup
			\end{table}

%%% INTRODUCE FIGURES AND DISCUSS RESULTS

% global variables
			\begin{figure}[t]
				\begin{adjustbox}{center}
					\begin{tabular}{cc}
						\begin{subfigure}{.5\textwidth}
  							\includegraphics[width=\linewidth]{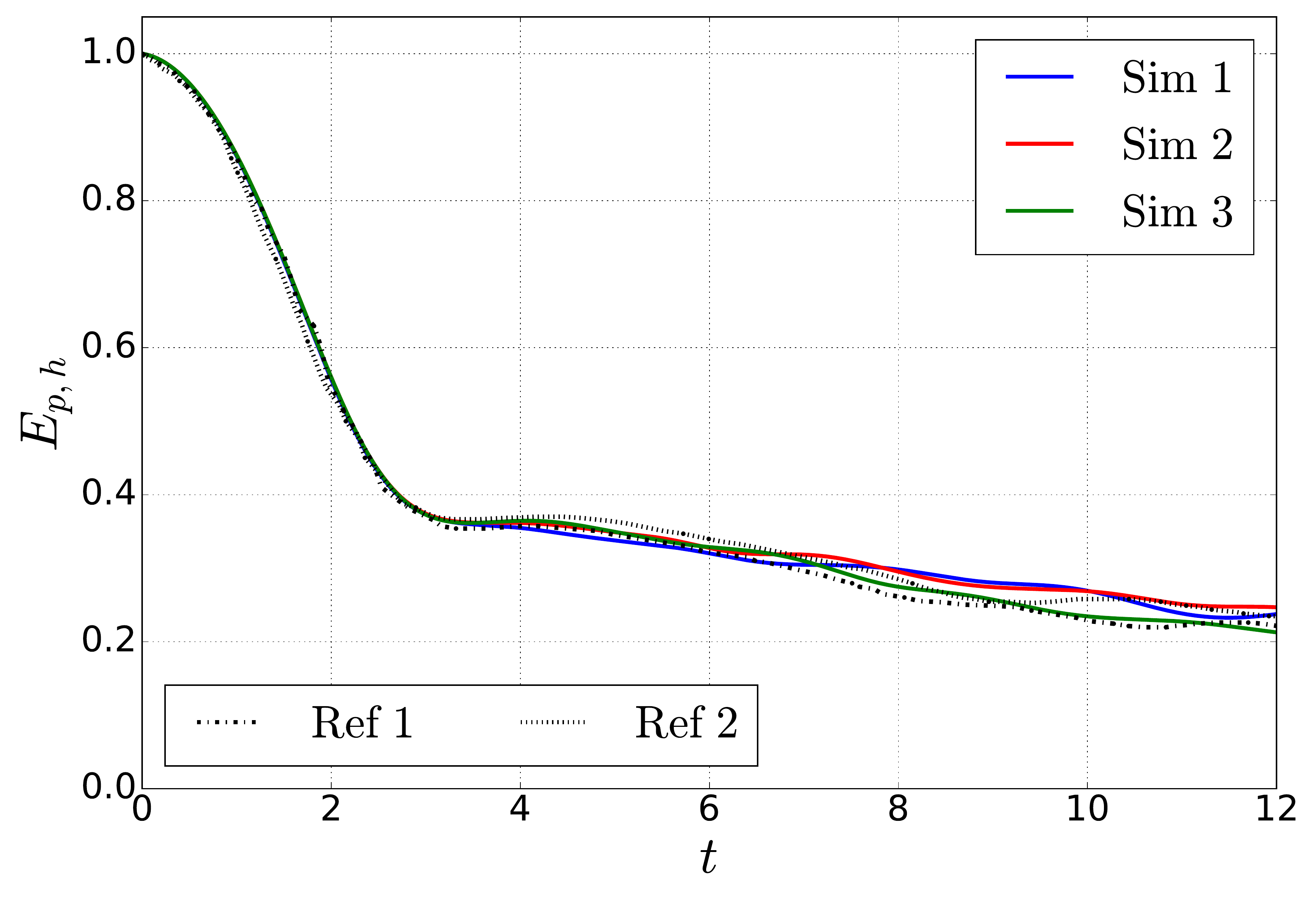}
							\caption[]{{\scriptsize Potential energy $E_{p,h}$.}}\label{subfig:ep}
						\end{subfigure}
						&
						\begin{subfigure}{.5\textwidth}
  							\includegraphics[width=\linewidth]{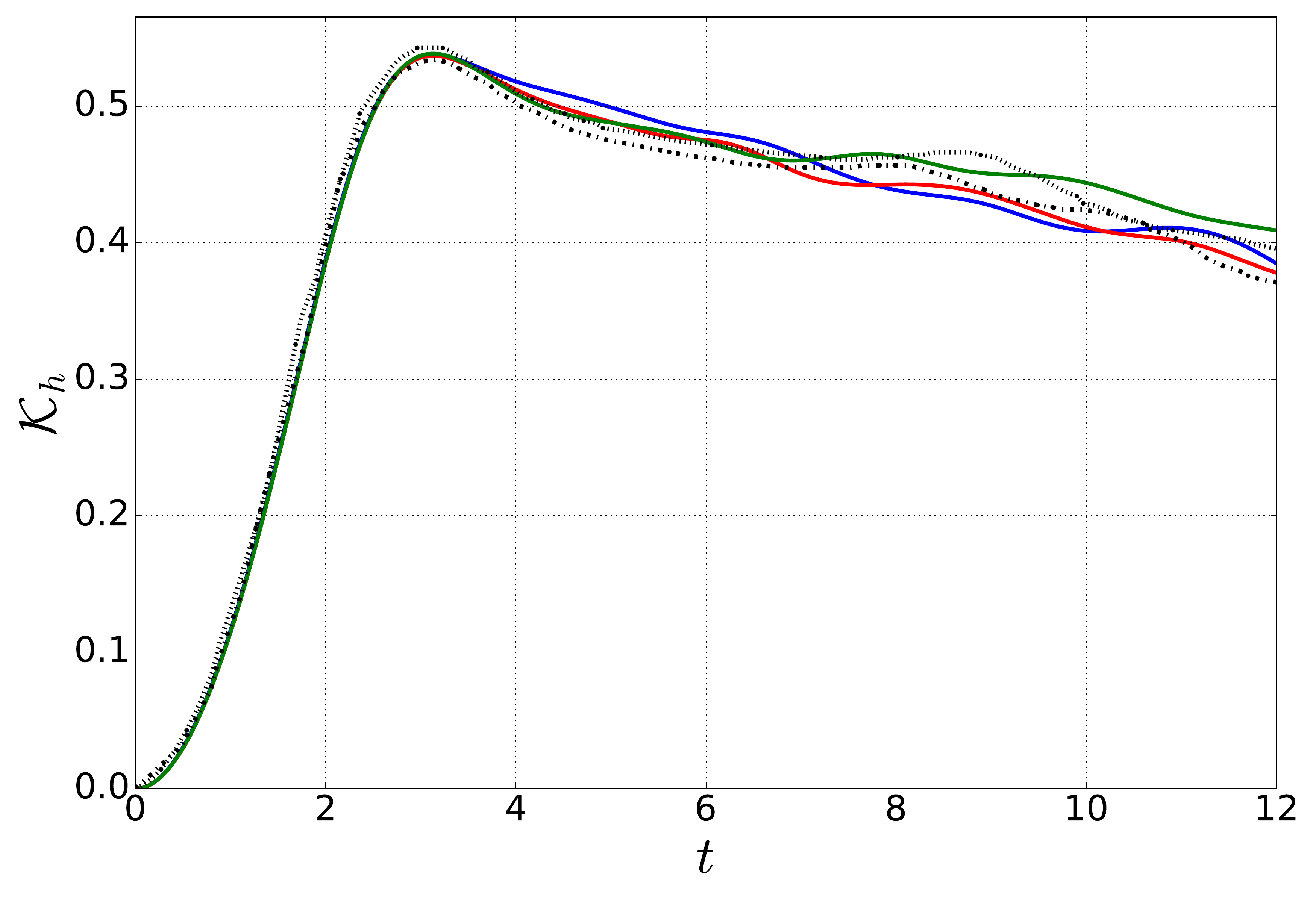}
							\caption[]{{\scriptsize Kinetic energy $\mathcal{K}_h$.}}\label{subfig:k}
						\end{subfigure}
						\\
						\begin{subfigure}{.5\textwidth}
  							\includegraphics[width=\linewidth]{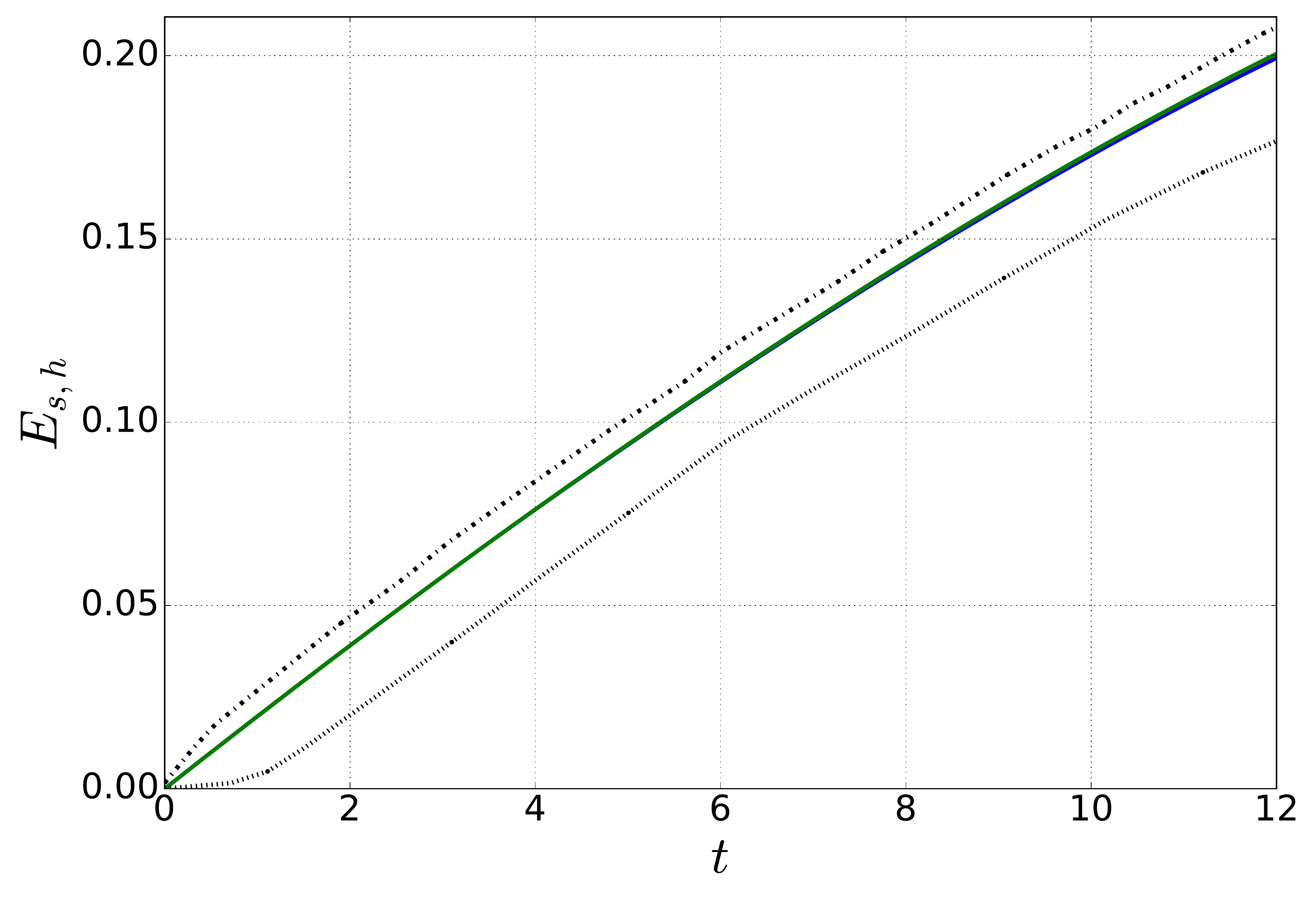}
							\caption[]{{\scriptsize Particle-settling dissipation $E_{s,h}$.}}\label{subfig:es}
						\end{subfigure}
						&
						\begin{subfigure}{.5\textwidth}
  							\includegraphics[width=\linewidth]{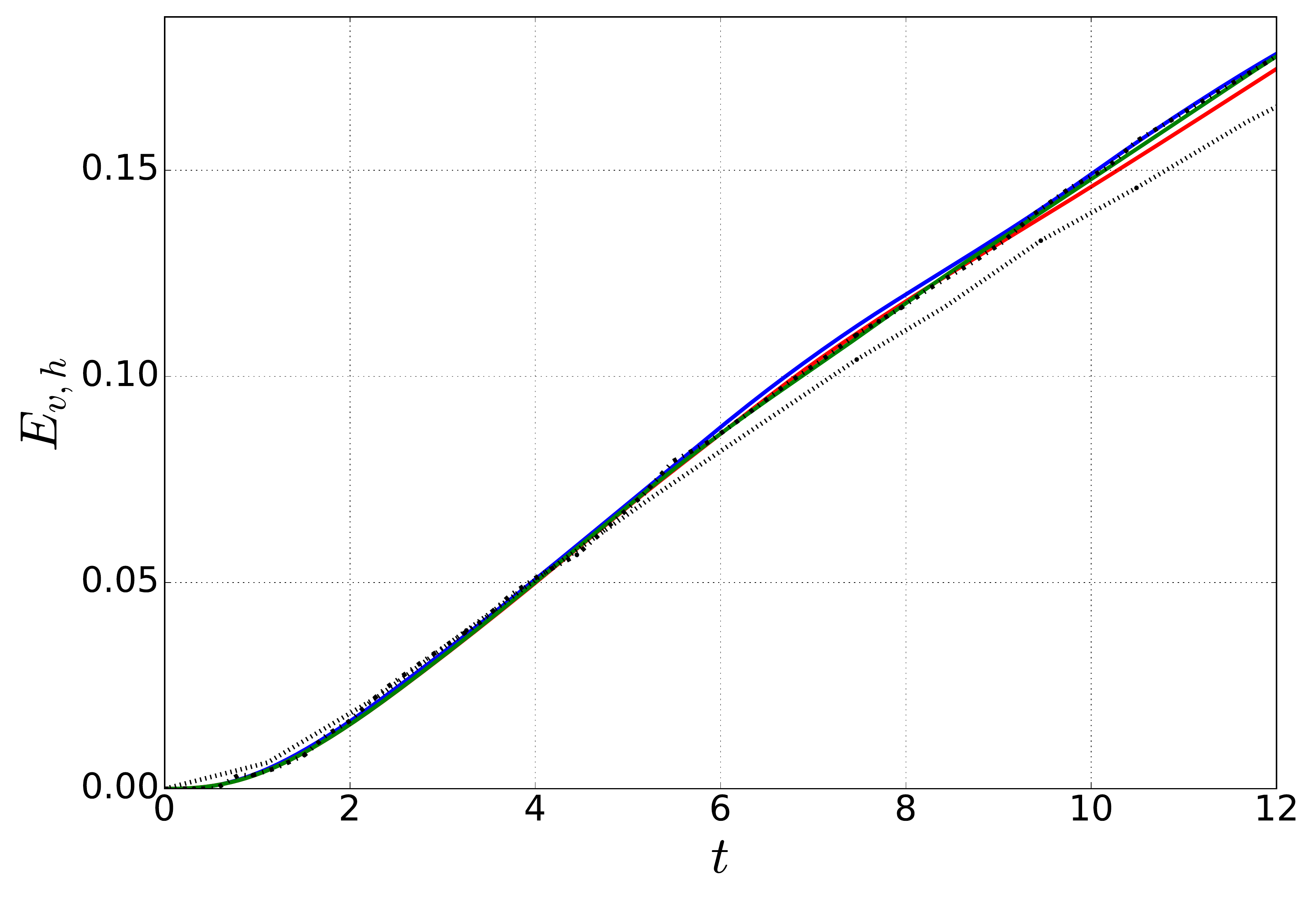}
							\caption[]{{\scriptsize Viscous dissipation $E_{v,h}$.}}\label{subfig:ev}
						\end{subfigure}
					\end{tabular}
				\end{adjustbox}
				\caption{Variation over time of the 4 components of the energy balance equation. Results are compared with those of \cite{parkinson2014,espath2014}.}
				\label{fig:energy_le}
			\end{figure}

% Comment on contours

			Figure \ref{fig:colorplot_le} displays the color plots of the particle concentration at $t=4$ and $12$, together with the corresponding meshes. Color plots from Parkinson et al. \cite{parkinson2014} are also shown because they illustrate the effects of mesh refinement and a mixed finite element method is also used. The first plot is a low resolution computation with a fixed mesh, while the second and third plots are high resolution computations with a fixed and an adaptive mesh, respectively. These last two plots display similar results, \cite{parkinson2014}.

			The particle concentration fields obtained with the modified MEEVC scheme resemble the high resolution reference results for $t=4$. For $t=4$, an intrusion front has fully developed and a vortical structure, characteristic of Kelvin-Helmholtz instabilities, is present along the upper interface. The reference results display two vortices along the upper interface, although for the low resolution case these are highly dissipated. With the MEEVC scheme these two vortices are captured, even with the lowest resolution. Interestingly, a third vortex appears in all of the MEEVC scheme computations. 

			For $t=12$, large variations appear from one computation to another. According to Parkinson et al. \cite{parkinson2014}, no two simulations ever produced the same results. This is due to the highly chaotic nature of the flow; vortices are generated and propagated in such a way that small variations in the mesh induce small variations in these vortices which grow and propagate over time. However, in the two highly-resolved computations from the reference, three distinct vortical structures, suspended over the remains of the intrusion front, can be observed. These structures can also be observed in the three computations carried out with the MEEVC scheme. Considerable variations in the shape and position of these structures are found from Sim 2 to Sim 3. 

% Energy
			In Figure \ref{fig:energy_le}, the evolution of the 4 discrete energy components can be found, see \eqref{eq:energy_eq_3}. The kinetic and potential energy are very similar for $t<3$. Until $t=6$, Sim 2 and 3 yield similar results. For $t>6$, vortex interaction becomes stronger and the three simulations differ considerably. It can be conlcuded that, in general, good agreement can be found between the MEEVC computations and those of the references. 

			For the particle-settling dissipation, hardly no difference is found between the three MEEVC simulations. The values for $E_{s,h}$ given by the references differ from each other and from the MEEVC scheme along the entire computation. Parkinson et al. use different boundary conditions at the top wall and the two references use a different expression for $\epsilon_{s}$. Although these are equivalent at the continuous level, the discrete expressions may differ. Parkinson et al.  calculate $\varepsilon_s$ with
			\begin{equation}\label{eq:vareps_s_park}
				\epsilon_{s,\mathrm{ref}\ 1} = - u_s \langle \vec{e}_g, \gradient \phi \rangle_{\Omega}- \frac{1}{\sqrt{\mathrm{Gr}\mathrm{Sc}^2}} \left\lbrace \langle \gradient \phi, \gradient y \rangle_{\Omega}- \int_{\partial \Omega} y \gradient \phi \cdot \vec{n}\ \mathrm{d}\Gamma \right\rbrace,
			\end{equation}

			\noindent and  Espath et al. use,
			\begin{equation}\label{eq:vareps_s_espath}
				\epsilon_{s,\mathrm{ref}\ 2} = - u_s \langle \vec{e}_g, \gradient \phi \rangle_{\Omega}+ \frac{1}{\sqrt{\mathrm{Gr}\mathrm{Sc}^2}} \langle \Delta \phi, y \rangle_\Omega.
			\end{equation}

\begin{figure}[]
				\centering
					\includegraphics[width=0.6\textwidth]{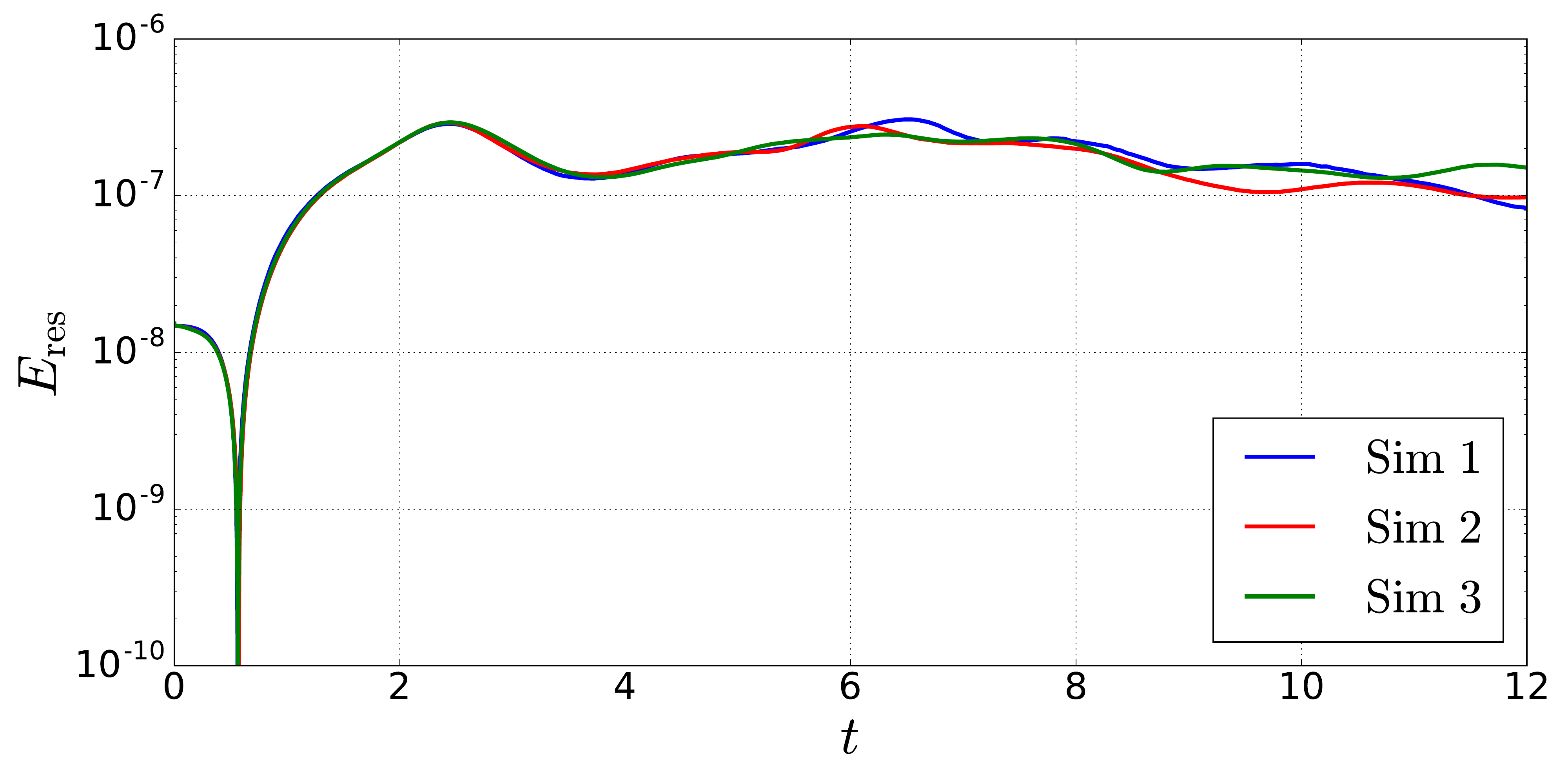}
					\caption{Energy residual $E_{\mathrm{res}}$.}
					\label{fig:res}
			\end{figure}

			A powerful advantage of the MEEVC scheme is the energy conserving properties. Figure \ref{fig:res} displays the evolution of the energy residual, $E_{res}$, defined,
			\begin{equation}\label{eq:energy_balance_5}
				E^k_{res} = \mathcal{K}_h^{k} + E_{p,h}^{k} +E_{v,h}^{k} + E^{k}_{s,h} - E_{p,h}^{0} - \mathcal{K}^0,
			\end{equation}

			\noindent where $\mathcal{K}_h^{k} = \frac{1}{2} ( \mathcal{K}_h^{k+\frac{3}{2}} + \mathcal{K}_h^{k+\frac{1}{2}})$. As expected, this term remains stable and bounded to relatively small values, indicating a lack of artificial dissipation. This error is only due to the staggering in time of the velocity, see \eqref{eq:energy_balance_4}. In contrast, Parkinson et al. \cite{parkinson2014} give the values resulting from integrating the energy errors over 60 time units, and these values are of the order of $10^{-1}-10^{-2}$, see Figure 3 from \cite{parkinson2014}.

\begin{figure}[]
				\begin{adjustbox}{center}
  					\begin{subfigure}{.35\textwidth}
    						\includegraphics[width=\linewidth]{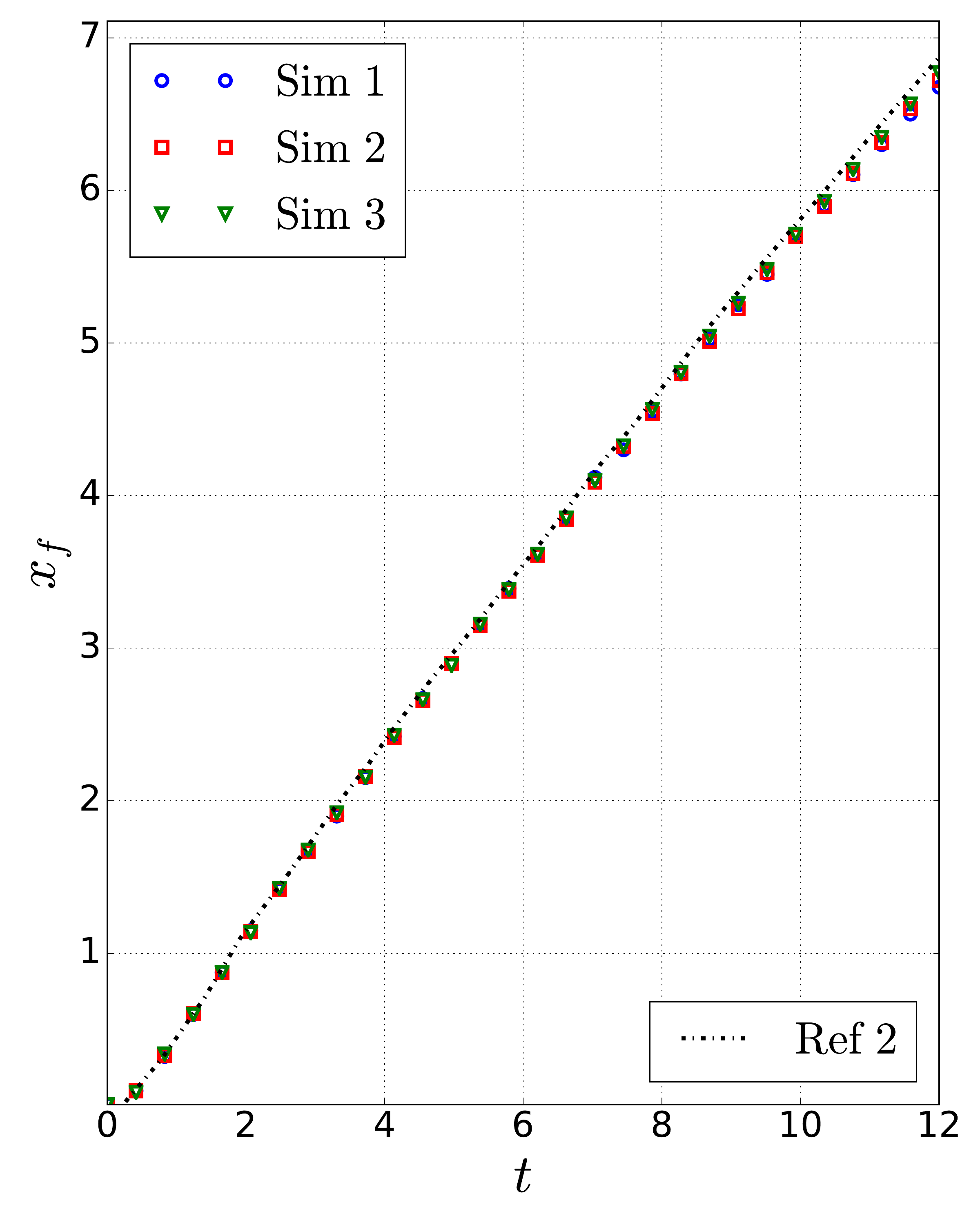}
  						\caption[]{{\scriptsize Position of the front $x_f$.}}     
    						\label{subfig:xf}
 	 				\end{subfigure}
  					\begin{subfigure}{.35\textwidth}
    						\includegraphics[width=\linewidth]{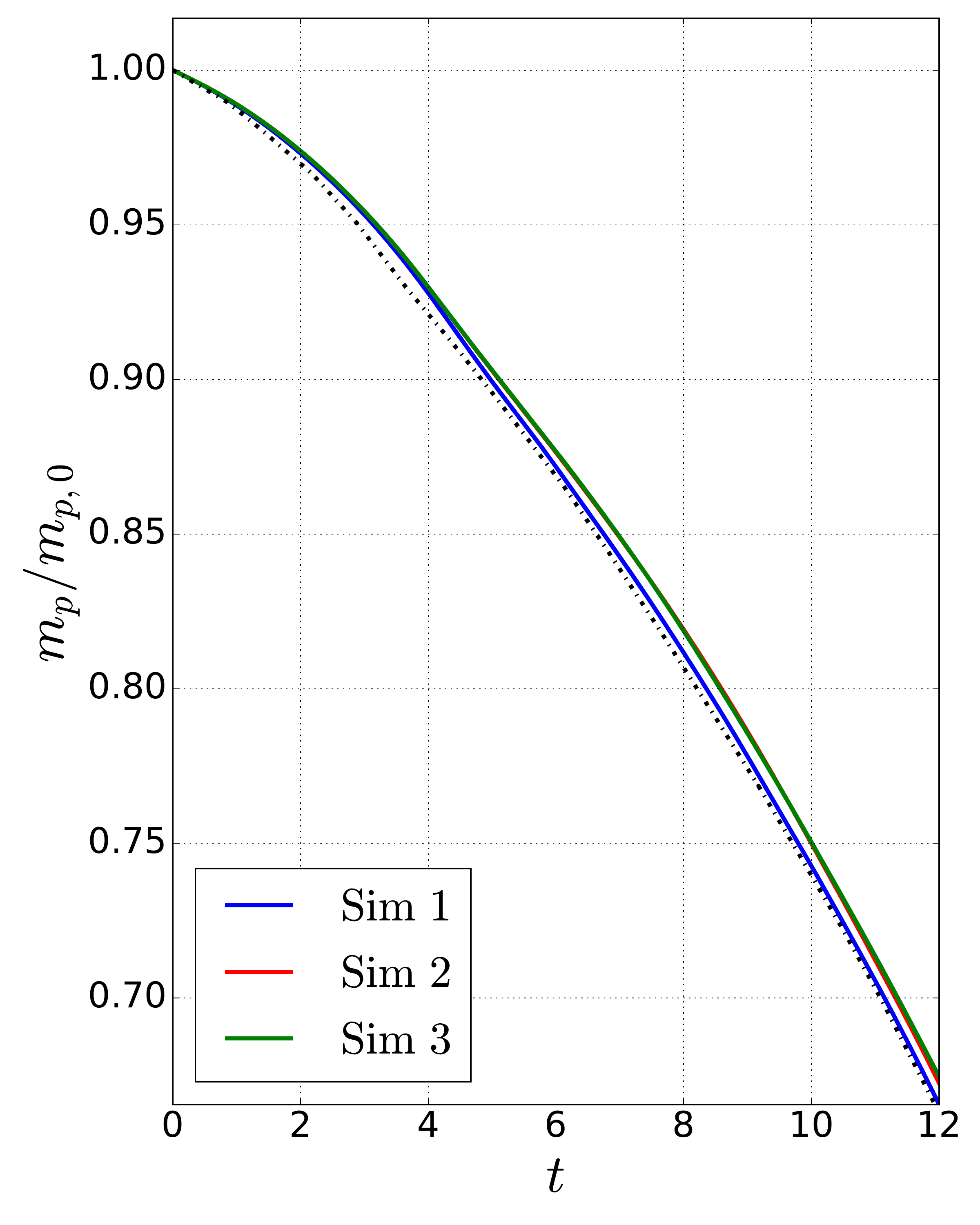}
  						\caption[]{{\scriptsize Variation of total mass $\frac{m_p}{m_{p,0}}$.}}      
    						\label{subfig:mp}
  					\end{subfigure}
  					\begin{subfigure}{.35\textwidth}
    						\includegraphics[width=\linewidth]{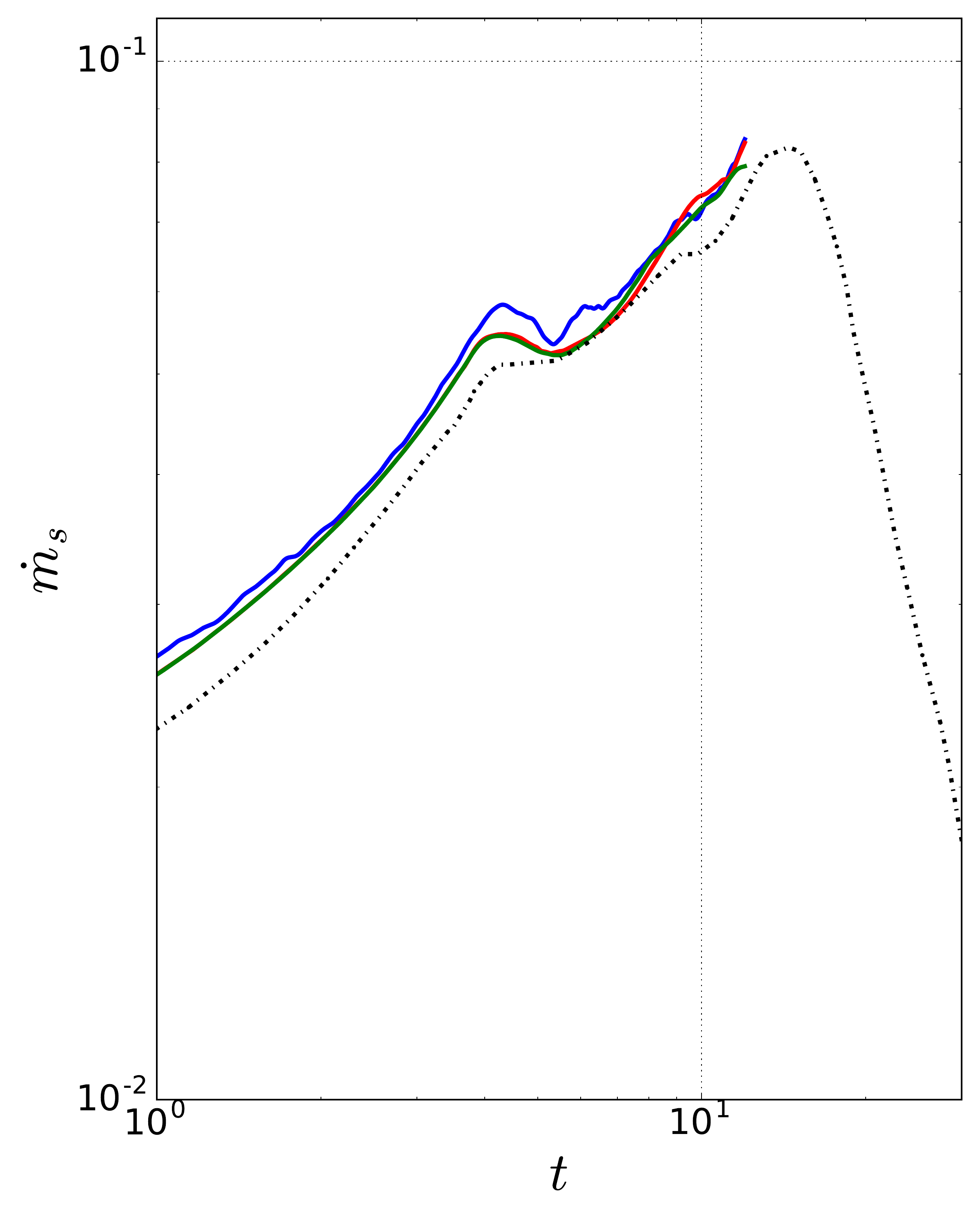}
      					\caption[]{{\scriptsize Sedimentation rate $\dot{m}_s$.}}      
    						\label{subfig:ms}
  					\end{subfigure}
  				\end{adjustbox}
  				\caption{Variation over time of the location of the front, the suspended mass and the sedimentation rate. Comparisons are made with the results from \citep{espath2014}.}
  				\label{fig:variables_2}
  			\end{figure}

\begin{sidewaysfigure}[]
				\begin{adjustbox}{center}\hspace{-1in}
					\begin{tabular}{cccc}
& \multicolumn{2}{c}{\scriptsize Results from Parkinson et al. \cite{parkinson2014}} & \\
	\scriptsize{\begin{tabular}{c} $ 2.9\times 10^4$ cells\\ (fixed mesh) \end{tabular}}
& \begin{subfigure}{.30\textwidth}\includegraphics[width=\linewidth]{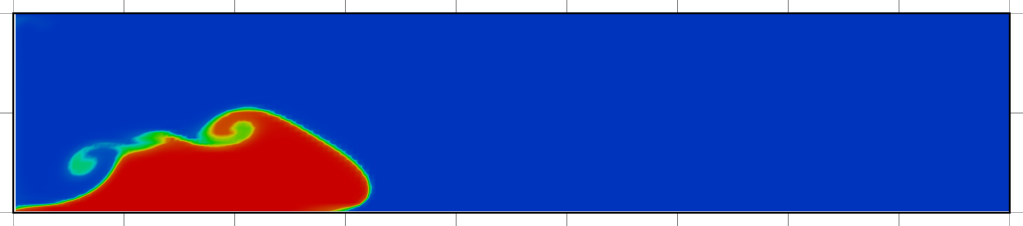}\end{subfigure} 
& \begin{subfigure}{.30\textwidth}\includegraphics[width=\linewidth]{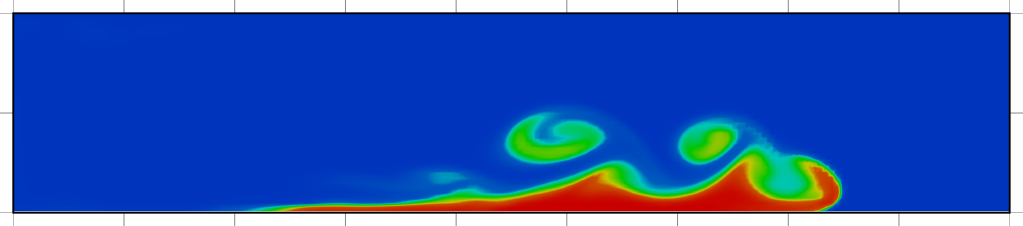}\end{subfigure} & \\ 
	\scriptsize{\begin{tabular}{c} $ 1.8\times 10^6$ cells\\ (fixed mesh) \end{tabular}}
& \begin{subfigure}{.30\textwidth}\includegraphics[width=\linewidth]{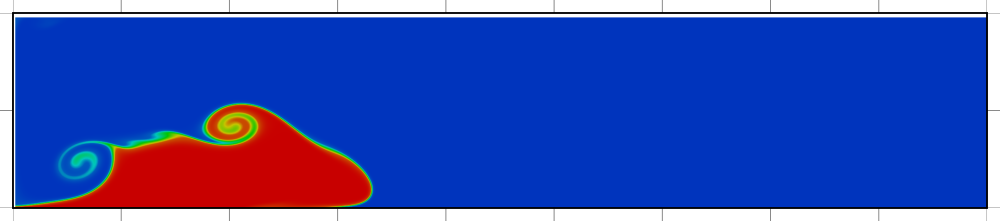}\end{subfigure}
& \begin{subfigure}{.30\textwidth}\includegraphics[width=\linewidth]{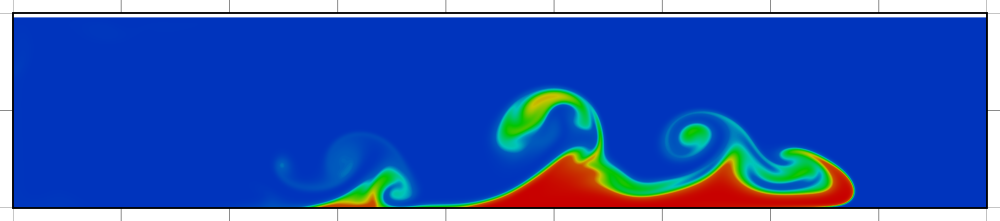}\end{subfigure} & \\ 
	\scriptsize{\begin{tabular}{c} $ 4.1\times 10^4$ cells\\ (adaptive mesh) \end{tabular}}
& \begin{subfigure}{.30\textwidth}\includegraphics[width=\linewidth]{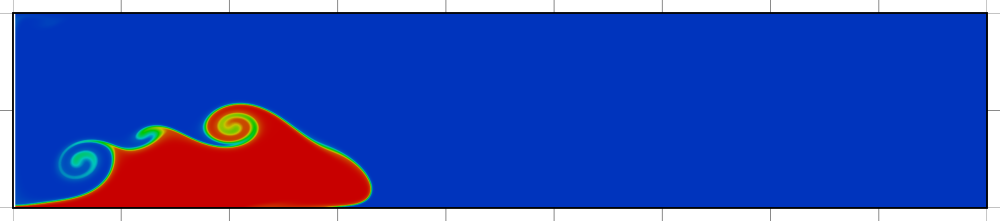}\end{subfigure}
& \begin{subfigure}{.30\textwidth}\includegraphics[width=\linewidth]{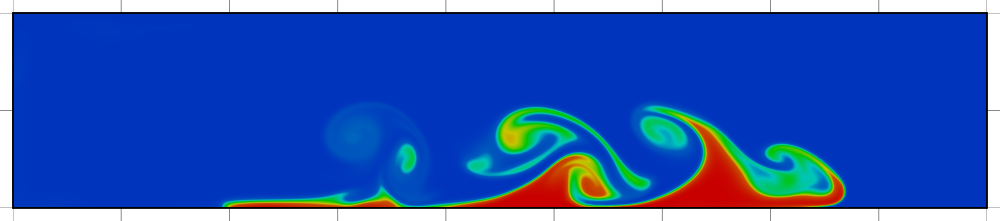}\end{subfigure} & \\ \\
& \multicolumn{2}{c}{\scriptsize Results obtained with the MEEVC scheme.} &  \\
	\scriptsize{\begin{tabular}{c} Sim 1 \\ $ 2.6\times 10^4$ cells (equivalent) \end{tabular}}
& \begin{subfigure}{.30\textwidth}\includegraphics[width=\linewidth]{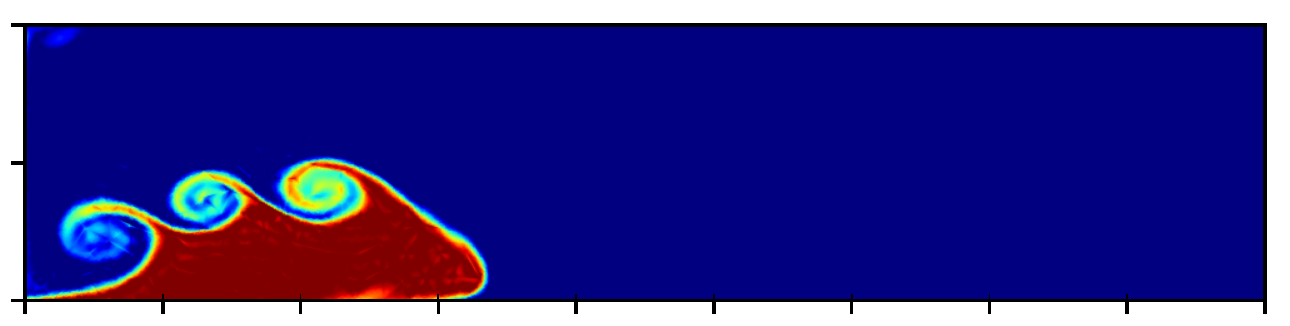}\end{subfigure}
& \begin{subfigure}{.30\textwidth}\includegraphics[width=\linewidth]{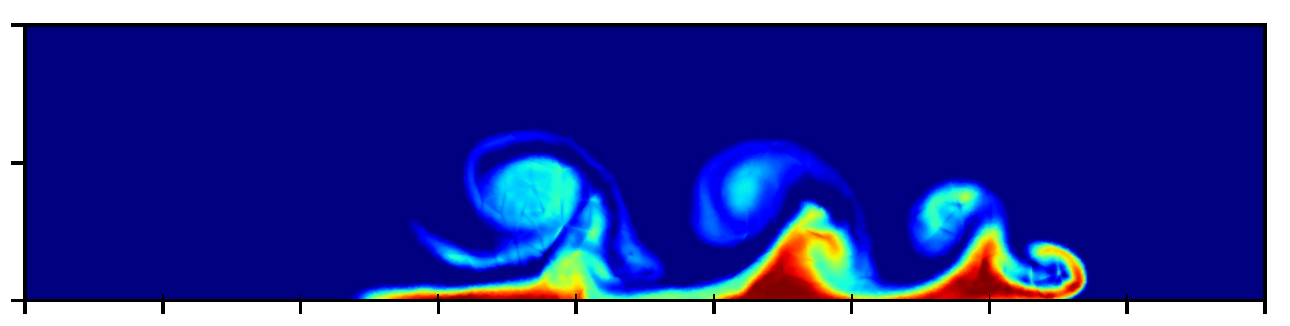}\end{subfigure}
& \begin{subfigure}{.30\textwidth}\includegraphics[width=\linewidth]{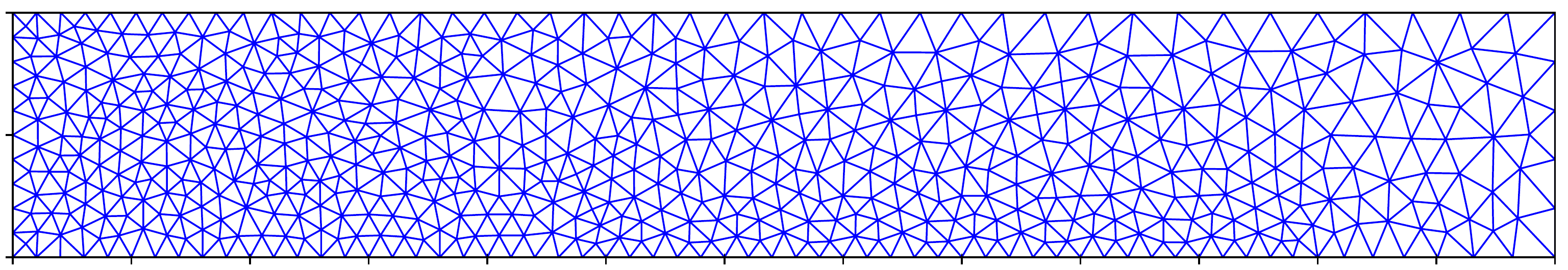}\end{subfigure} \\
	\scriptsize{\begin{tabular}{c} Sim 2 \\ $ 5.6\times 10^4$ cells (equivalent) \end{tabular}}
& \begin{subfigure}{.30\textwidth}\includegraphics[width=\linewidth]{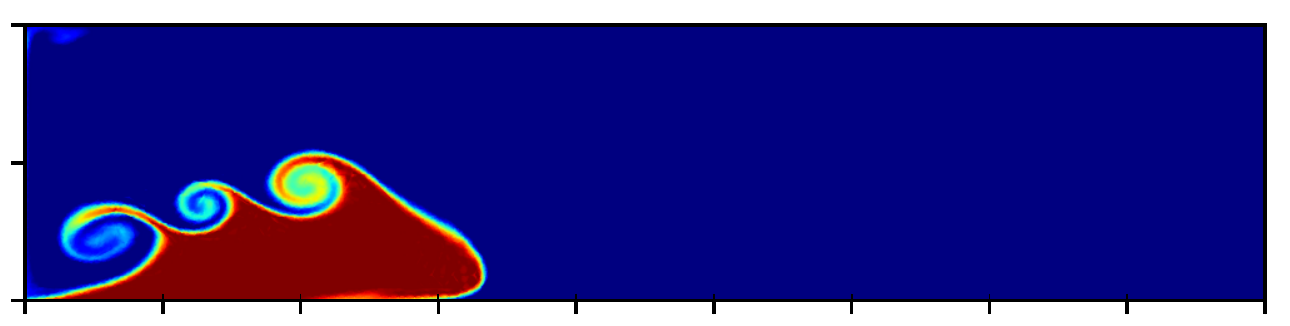}\end{subfigure}
& \begin{subfigure}{.30\textwidth}\includegraphics[width=\linewidth]{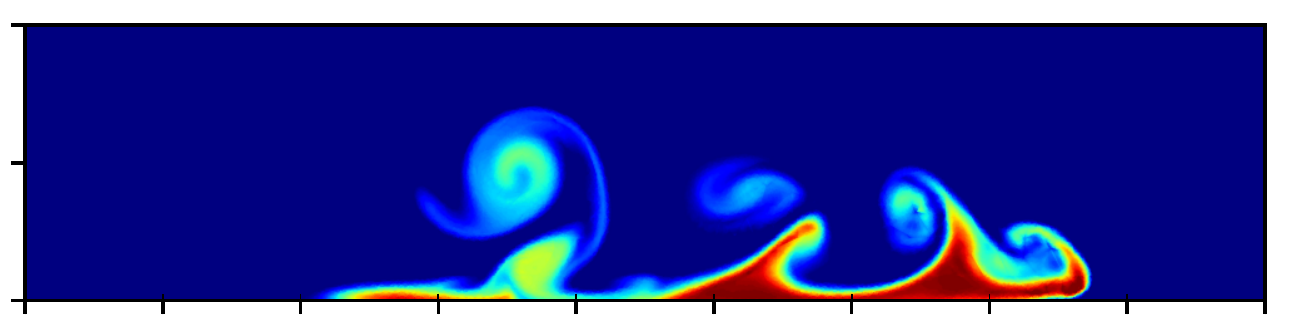}\end{subfigure}
& \begin{subfigure}{.30\textwidth}\includegraphics[width=\linewidth]{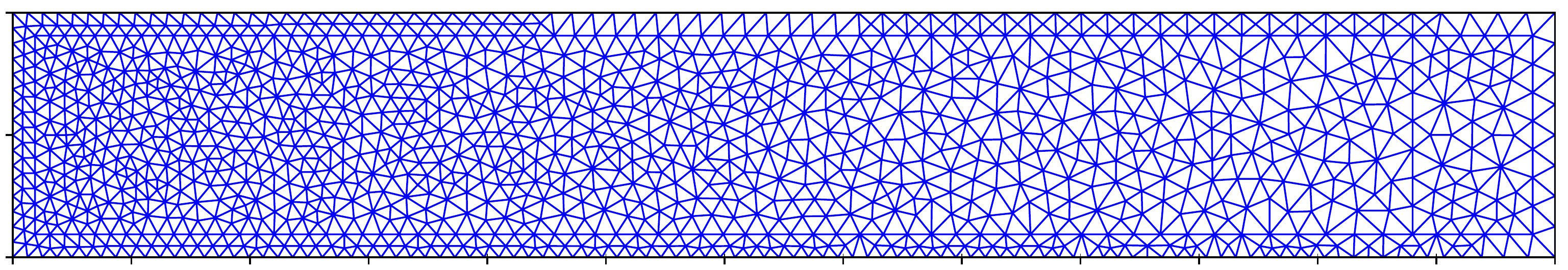}\end{subfigure} \\
	\scriptsize{\begin{tabular}{c} Sim 3 \\ $ 8.4 \times 10^4$ cells (equivalent) \end{tabular}}
& \begin{subfigure}{.30\textwidth}\includegraphics[width=\linewidth]{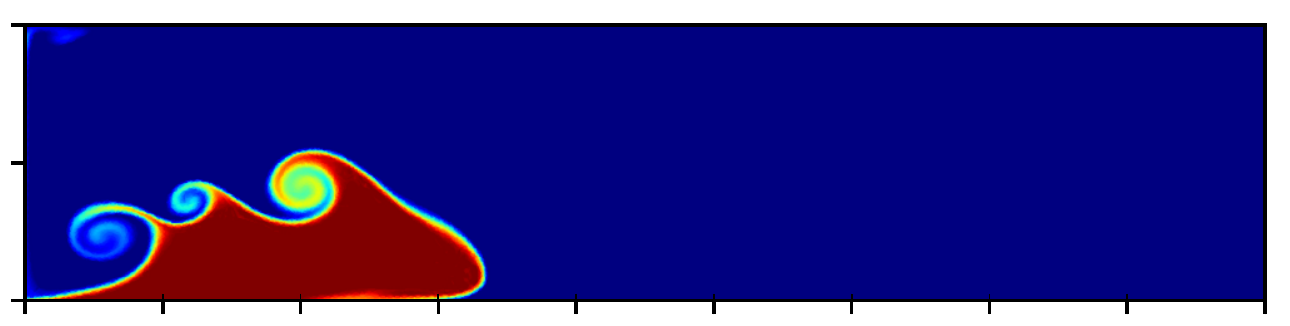}\end{subfigure}
& \begin{subfigure}{.30\textwidth}\includegraphics[width=\linewidth]{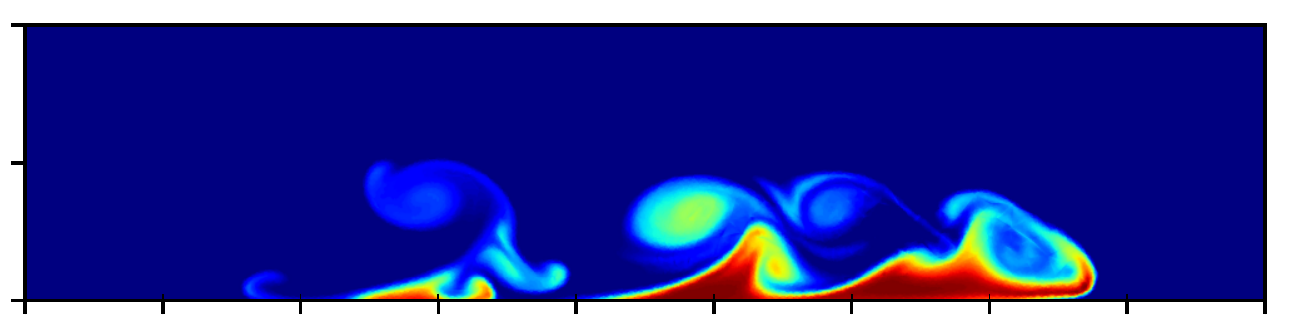}\end{subfigure}
& \begin{subfigure}{.30\textwidth}\includegraphics[width=\linewidth]{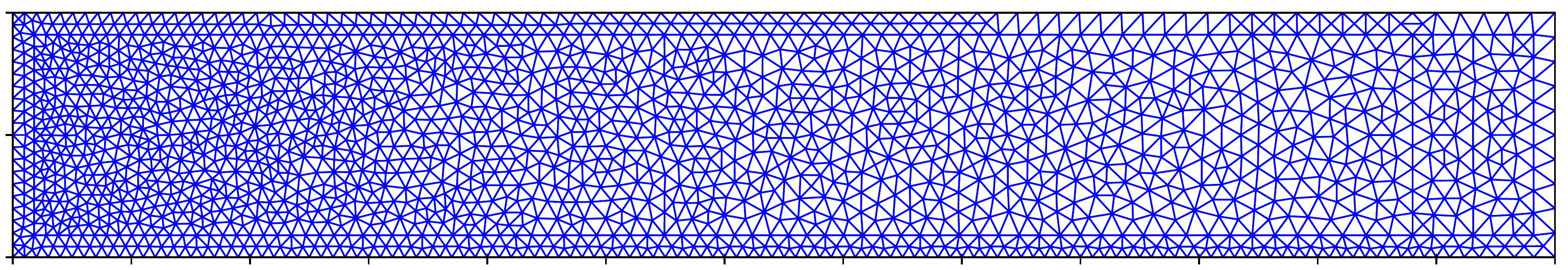}\end{subfigure} \\
& \scriptsize $t = 4$ & \scriptsize $t = 12$ & \\
& \multicolumn{2}{c}{ \includegraphics[width=0.5\linewidth]{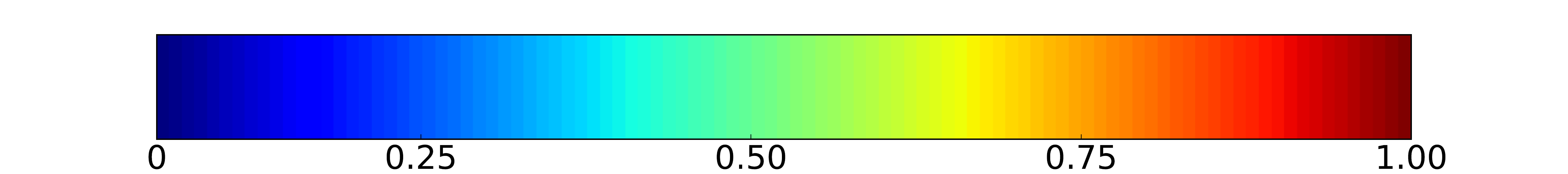} } &  \\
\end{tabular}
\end{adjustbox}
\caption{Distribution of particle concentration at $t=4$ and $12$ obtained with the three discretizations described in Table \ref{tab:le_sim}, together with the corresponding meshes. Three sets of results from \cite{parkinson2014} are also presented.}
\label{fig:colorplot_le}
\end{sidewaysfigure}

			Three of the main features of a turbidity current are the location of the front $x_f$, the mass of suspended particles $m_p$ and the sedimentation rate $\dot{m}_s$. The evolution of these variables is shown in Figure \ref{fig:variables_2}, together with data from Espath et al. \cite{espath2014}. The calculated values for the location of the front and the total suspended mass agree with the reference results. This last variable is defined as
			\begin{equation}\label{eq:mp}
				\frac{m_p}{m_{p,0}} = \frac{1}{m_{p,0}} \int_\Omega \phi\ \mathrm{d}\Omega, \quad \text{with} \quad m_{p,0} = \int_\Omega \phi(t=0)\ \mathrm{d}\Omega.
			\end{equation}

			The variation of $m_p$, denoted $\dot{m}_s$, is defined as
			\begin{equation}\label{eq:ms}
				\dot{m}_s = - \int_{\Gamma_3} \phi u_s\ \mathrm{d}\Gamma,
			\end{equation}

			\noindent see \eqref{subsec:le_bc} for a derivation. Figure \ref{subfig:ms} indicates that the results from Sim 2 and 3 are in great agreement and that a similar tendency to that of the reference is obtained. On the other hand, the results from Sim 1 show a more chaotic behavior.

% global variables 2

% global variables

	\section{Concluding remarks}

		In this document, a physics-compatible solver for turbidity currents is constructed from the MEEVC scheme with kinematic Neumann boundary conditions. It is proved that the discrete energy balance equation holds up to a residual due to the staggering in time of the velocity which does not accumulate over time and is proportional to the time step $\Delta t$. With the finite element basis functions used for the velocity and the pressure, $\divergence \vec{u}_h = 0$ in $\Omega$ and certain vector calculus identities required for the derivation of the energy balance equation hold at the discrete level. Essentially, the discrete kinetic and potential energy are only lost to the viscous and sedimentation terms, and no artificial dissipation of energy is introduced. 

		This solver is validated by computing the lock exchange test case for $\mathrm{Gr} = 5\times10^6$ ($\mathrm{Re} \approx 2236$) with three different meshes of increasing refinement. Comparisons are made with Parkinson et al. \cite{parkinson2014} and Espath et al. \cite{espath2014}, and these indicate that the essential dynamics and vortical structures are captured with the lowest resolution. Further features of the flow are well computed and compared with the reference results, showing very good agreement in general.
		
		The dual field formulation is currently being extended to 3D flows.

%% References with bibTeX database -------------------------------------------------------------------
	\bibliographystyle{model1-num-names}
	\bibliography{Bibliography.bib}

\begin{thebibliography}{49}
\expandafter\ifx\csname natexlab\endcsname\relax\def\natexlab#1{#1}\fi
\providecommand{\bibinfo}[2]{#2}
\ifx\xfnm\relax \def\xfnm[#1]{\unskip,\space#1}\fi
%Type = Article
\bibitem[{Palha and Gerritsma(2017)}]{palha2016}
\bibinfo{author}{A.~Palha}, \bibinfo{author}{M.~I. Gerritsma},
\newblock \bibinfo{title}{A mass, energy, enstrophy and vorticity conserving
  ({MEEVC}) mimetic spectral element discretization for the {2D} incompressible
  {Navier-Stokes} equations},
\newblock \bibinfo{journal}{Journal of Computational Physics}
  \bibinfo{volume}{328} (\bibinfo{year}{2017}) \bibinfo{pages}{200--220}.
%Type = Article
\bibitem[{Simpson(1982)}]{Simpson1982}
\bibinfo{author}{J.~E. Simpson},
\newblock \bibinfo{title}{{Gravity Currents in the Laboratory, Atmosphere, and
  Ocean}},
\newblock \bibinfo{journal}{Annual Review of Fluid Mechanics}
  \bibinfo{volume}{14} (\bibinfo{year}{1982}) \bibinfo{pages}{213--234}.
%Type = Book
\bibitem[{Ungarish(2009)}]{Ungarish2009}
\bibinfo{author}{M.~Ungarish}, \bibinfo{title}{{An introduction to gravity
  currents and intrusions}}, \bibinfo{publisher}{CRC Press},
  \bibinfo{year}{2009}.
%Type = Article
\bibitem[{Meiburg et~al.(2015)Meiburg, Radhakrishnan, and
  Nasr-Azadani}]{Meiburg2015}
\bibinfo{author}{E.~Meiburg}, \bibinfo{author}{S.~Radhakrishnan},
  \bibinfo{author}{M.~Nasr-Azadani},
\newblock \bibinfo{title}{Modeling gravity and turbidity currents:
  computational approaches and challenges},
\newblock \bibinfo{journal}{Applied Mechanics Reviews} \bibinfo{volume}{67}
  (\bibinfo{year}{2015}) \bibinfo{pages}{040802:1--23}.
%Type = Article
\bibitem[{Huppert(2006)}]{Huppert2006299}
\bibinfo{author}{H.~Huppert},
\newblock \bibinfo{title}{Gravity currents: A personal perspective},
\newblock \bibinfo{journal}{Journal of Fluid Mechanics} \bibinfo{volume}{554}
  (\bibinfo{year}{2006}) \bibinfo{pages}{299--322}.
%Type = Article
\bibitem[{Meiburg and Kneller(2010)}]{Meiburg2010}
\bibinfo{author}{E.~Meiburg}, \bibinfo{author}{B.~Kneller},
\newblock \bibinfo{title}{{Turbidity Currents and Their Deposits}},
\newblock \bibinfo{journal}{Annual Review of Fluid Mechanics}
  \bibinfo{volume}{42} (\bibinfo{year}{2010}) \bibinfo{pages}{135--156}.
%Type = Article
\bibitem[{Nasr-Azadani and Meiburg(2014)}]{Nasr-Azadani2014}
\bibinfo{author}{M.~M. Nasr-Azadani}, \bibinfo{author}{E.~Meiburg},
\newblock \bibinfo{title}{{Turbidity currents interacting with
  three-dimensional seafloor topography}},
\newblock \bibinfo{journal}{Journal of Fluid Mechanics} \bibinfo{volume}{745}
  (\bibinfo{year}{2014}) \bibinfo{pages}{409--443}.
%Type = Article
\bibitem[{Kneller and McCaffrey(1999)}]{Kneller1999}
\bibinfo{author}{B.~Kneller}, \bibinfo{author}{W.~McCaffrey},
\newblock \bibinfo{title}{{Depositional effects of flow nonuniformity and
  stratification within turbidity currents approaching a bounding slope;
  deflection, reflection, and facies variation}},
\newblock \bibinfo{journal}{Journal of Sedimentary Research}
  \bibinfo{volume}{69} (\bibinfo{year}{1999}) \bibinfo{pages}{980--991}.
%Type = Article
\bibitem[{Syvitski et~al.(1996)Syvitski, Alexander, Field, Gardner, Orange, and
  Yun}]{Syvitski1996}
\bibinfo{author}{J.~Syvitski}, \bibinfo{author}{C.~Alexander},
  \bibinfo{author}{M.~Field}, \bibinfo{author}{J.~Gardner},
  \bibinfo{author}{D.~Orange}, \bibinfo{author}{J.~Yun},
\newblock \bibinfo{title}{{Continental-Slope Sedimentation: The View from
  Northern California}},
\newblock \bibinfo{journal}{Oceanography} \bibinfo{volume}{9}
  (\bibinfo{year}{1996}) \bibinfo{pages}{163--167}.
%Type = Article
\bibitem[{Parkinson et~al.(2014)Parkinson, Hill, Piggott, and
  Allison}]{parkinson2014}
\bibinfo{author}{S.~D. Parkinson}, \bibinfo{author}{J.~Hill},
  \bibinfo{author}{M.~D. Piggott}, \bibinfo{author}{P.~A. Allison},
\newblock \bibinfo{title}{Direct numerical simulations of particle-laden
  density currents with adaptive, discontinuous finite elements},
\newblock \bibinfo{journal}{Geoscientific Model Development}
  \bibinfo{volume}{7} (\bibinfo{year}{2014}) \bibinfo{pages}{1945--1960}.
%Type = Article
\bibitem[{Kneller and Buckee(2000)}]{Kneller2000}
\bibinfo{author}{B.~Kneller}, \bibinfo{author}{C.~Buckee},
\newblock \bibinfo{title}{{The structure and fluid mechanics of turbidity
  currents: a review of some recent studies and their geological
  implications}},
\newblock \bibinfo{journal}{Sedimentology} \bibinfo{volume}{47}
  (\bibinfo{year}{2000}) \bibinfo{pages}{62--94}.
%Type = Article
\bibitem[{von K\'{a}rm\'{a}n(1940)}]{vonKarman1940}
\bibinfo{author}{T.~von K\'{a}rm\'{a}n},
\newblock \bibinfo{title}{{The engineer grapples with non-linear problems}},
\newblock \bibinfo{journal}{Bulletin American Mathematical Society}
  \bibinfo{volume}{46} (\bibinfo{year}{1940}) \bibinfo{pages}{615--683}.
%Type = Article
\bibitem[{Konopliv et~al.(2016)Konopliv, {Llewellyn Smith}, McElwaine, and
  Meiburg}]{Konopliv2016}
\bibinfo{author}{N.~A. Konopliv}, \bibinfo{author}{S.~G. {Llewellyn Smith}},
  \bibinfo{author}{J.~N. McElwaine}, \bibinfo{author}{E.~Meiburg},
\newblock \bibinfo{title}{{Modelling gravity currents without an
  energy closure}},
\newblock \bibinfo{journal}{Journal of Fluid Mechanics} \bibinfo{volume}{789}
  (\bibinfo{year}{2016}) \bibinfo{pages}{806--829}.
%Type = Book
\bibitem[{Kundu et~al.(2012)Kundu, Cohen, and Dowling}]{Kundu2012}
\bibinfo{author}{P.~K. Kundu}, \bibinfo{author}{I.~M. Cohen},
  \bibinfo{author}{D.~R. Dowling}, \bibinfo{title}{{Fluid mechanics}},
  \bibinfo{publisher}{Academic Press}, \bibinfo{year}{2012}.
%Type = Article
\bibitem[{Spalart(2000)}]{Spalart2000}
\bibinfo{author}{P.~R. Spalart},
\newblock \bibinfo{title}{{Strategies for turbulence modelling and
  simulations}},
\newblock \bibinfo{journal}{International Journal of Heat and Fluid Flow}
  \bibinfo{volume}{21} (\bibinfo{year}{2000}) \bibinfo{pages}{252--263}.
%Type = Article
\bibitem[{Lesieur and Metais(1996)}]{Lesieur1996}
\bibinfo{author}{M.~Lesieur}, \bibinfo{author}{O.~Metais},
\newblock \bibinfo{title}{{New Trends in Large-Eddy Simulations of
  Turbulence}},
\newblock \bibinfo{journal}{Annual Review of Fluid Mechanics}
  \bibinfo{volume}{28} (\bibinfo{year}{1996}) \bibinfo{pages}{45--82}.
%Type = Article
\bibitem[{Zhiyin(2015)}]{Zhiyin2015}
\bibinfo{author}{Y.~Zhiyin},
\newblock \bibinfo{title}{{Large-eddy simulation: Past, present and the
  future}},
\newblock \bibinfo{journal}{Chinese Journal of Aeronautics}
  \bibinfo{volume}{28} (\bibinfo{year}{2015}) \bibinfo{pages}{11--24}.
%Type = Article
\bibitem[{H{\"a}rtel et~al.(1997)H{\"a}rtel, Kleiser, Michaud, and
  Stein}]{hartel1997}
\bibinfo{author}{C.~H{\"a}rtel}, \bibinfo{author}{L.~Kleiser},
  \bibinfo{author}{M.~Michaud}, \bibinfo{author}{C.~F. Stein},
\newblock \bibinfo{title}{A direct numerical simulation approach to the study
  of intrusion flows},
\newblock \bibinfo{journal}{Journal of Engineering Mathematics}
  \bibinfo{volume}{32} (\bibinfo{year}{1997}) \bibinfo{pages}{103--120}.
%Type = Article
\bibitem[{H{\"a}rtel et~al.(2000)H{\"a}rtel, Meiburg, and Necker}]{hartel2000}
\bibinfo{author}{C.~H{\"a}rtel}, \bibinfo{author}{E.~Meiburg},
  \bibinfo{author}{F.~Necker},
\newblock \bibinfo{title}{Analysis and direct numerical simulation of the flow
  at a gravity-current head. part 1. flow topology and front speed for slip and
  no-slip boundaries},
\newblock \bibinfo{journal}{Journal of Fluid Mechanics} \bibinfo{volume}{418}
  (\bibinfo{year}{2000}) \bibinfo{pages}{189--212}.
%Type = Article
\bibitem[{Necker et~al.(2002)Necker, H{\"a}rtel, Kleiser, and
  Meiburg}]{necker2002}
\bibinfo{author}{F.~Necker}, \bibinfo{author}{C.~H{\"a}rtel},
  \bibinfo{author}{L.~Kleiser}, \bibinfo{author}{E.~Meiburg},
\newblock \bibinfo{title}{High-resolution simulations of particle-driven
  gravity currents},
\newblock \bibinfo{journal}{International Journal of Multiphase Flow}
  \bibinfo{volume}{28} (\bibinfo{year}{2002}) \bibinfo{pages}{279--300}.
%Type = Article
\bibitem[{Necker et~al.(2005)Necker, H{\"a}rtel, Kleiser, and
  Meiburg}]{necker2005}
\bibinfo{author}{F.~Necker}, \bibinfo{author}{C.~H{\"a}rtel},
  \bibinfo{author}{L.~Kleiser}, \bibinfo{author}{E.~Meiburg},
\newblock \bibinfo{title}{Mixing and dissipation in particle-driven gravity
  currents},
\newblock \bibinfo{journal}{Journal of Fluid Mechanics} \bibinfo{volume}{545}
  (\bibinfo{year}{2005}) \bibinfo{pages}{339--372}.
%Type = Article
\bibitem[{Blanchette et~al.(2005)Blanchette, Strauss, Meiburg, Kneller, and
  Glinsky}]{blanchette2005}
\bibinfo{author}{F.~Blanchette}, \bibinfo{author}{M.~Strauss},
  \bibinfo{author}{E.~Meiburg}, \bibinfo{author}{B.~Kneller},
  \bibinfo{author}{M.~E. Glinsky},
\newblock \bibinfo{title}{High-resolution numerical simulations of resuspending
  gravity currents: Conditions for self-sustainment},
\newblock \bibinfo{journal}{Journal of Geophysical Research}
  \bibinfo{volume}{110} (\bibinfo{year}{2005}) \bibinfo{pages}{C12022}.
%Type = Article
\bibitem[{Cantero et~al.(2006)Cantero, Balachandar, Garcia, and
  Ferry}]{cantero2006}
\bibinfo{author}{M.~I. Cantero}, \bibinfo{author}{S.~Balachandar},
  \bibinfo{author}{M.~H. Garcia}, \bibinfo{author}{J.~P. Ferry},
\newblock \bibinfo{title}{Direct numerical simulations of planar and
  cylindrical density currents},
\newblock \bibinfo{journal}{Journal of Applied Mechanics} \bibinfo{volume}{73}
  (\bibinfo{year}{2006}) \bibinfo{pages}{923--930}.
%Type = Article
\bibitem[{Birman and Meiburg(2005)}]{birman2005}
\bibinfo{author}{J.~E. Birman, V.\ K.~Martin}, \bibinfo{author}{E.~Meiburg},
\newblock \bibinfo{title}{The non-boussinesq lock-exchange problem. part 2.
  high-resolution simulations},
\newblock \bibinfo{journal}{Journal of Fluid Mechanics} \bibinfo{volume}{537}
  (\bibinfo{year}{2005}) \bibinfo{pages}{125--144}.
%Type = Article
\bibitem[{Cantero et~al.(2008)Cantero, Balachandar, and Garcia}]{cantero2008}
\bibinfo{author}{M.~I. Cantero}, \bibinfo{author}{S.~Balachandar},
  \bibinfo{author}{M.~H. Garcia},
\newblock \bibinfo{title}{An eulerian-eulerian model for gravity currents
  driven by inertial particles},
\newblock \bibinfo{journal}{International Journal of Multiphase Flow}
  \bibinfo{volume}{34} (\bibinfo{year}{2008}) \bibinfo{pages}{484--501}.
%Type = Article
\bibitem[{Espath et~al.(2014)Espath, Pinto, Laizet, and
  Silvestrini}]{espath2014}
\bibinfo{author}{L.~F.~R. Espath}, \bibinfo{author}{L.~C. Pinto},
  \bibinfo{author}{S.~Laizet}, \bibinfo{author}{J.~H. Silvestrini},
\newblock \bibinfo{title}{Two- and three-dimensional direct numerical
  simulation of particle-laden gravity currents},
\newblock \bibinfo{journal}{Computers and Geosciences} \bibinfo{volume}{63}
  (\bibinfo{year}{2014}) \bibinfo{pages}{9--16}.
%Type = Article
\bibitem[{de~Diego et~al.(2019)de~Diego, Palha, and Gerritsma}]{DeDiego2019a}
\bibinfo{author}{G.~G. de~Diego}, \bibinfo{author}{A.~Palha},
  \bibinfo{author}{M.~Gerritsma},
\newblock \bibinfo{title}{{Inclusion of no-slip boundary conditions in the
  MEEVC scheme}},
\newblock \bibinfo{journal}{Journal of Computational Physics}
  \bibinfo{volume}{378} (\bibinfo{year}{2019}) \bibinfo{pages}{615--633}.
%Type = Article
\bibitem[{Ferry and Balachandar(2001)}]{ferry2001}
\bibinfo{author}{J.~Ferry}, \bibinfo{author}{S.~Balachandar},
\newblock \bibinfo{title}{A fast {Eulerian} method for disperse two-phase
  flow},
\newblock \bibinfo{journal}{International Journal of Multiphase Flow}
  \bibinfo{volume}{27} (\bibinfo{year}{2001}) \bibinfo{pages}{1199--1226}.
%Type = Article
\bibitem[{Elghobashi(1994)}]{elghobashi1994}
\bibinfo{author}{S.~Elghobashi},
\newblock \bibinfo{title}{On predicting particle-laden turbulent flows},
\newblock \bibinfo{journal}{Applied Scientific Research} \bibinfo{volume}{52}
  (\bibinfo{year}{1994}) \bibinfo{pages}{309--329}.
%Type = Article
\bibitem[{Balachandar(2009)}]{balachandar2009}
\bibinfo{author}{S.~Balachandar},
\newblock \bibinfo{title}{A scaling analysis for point-particle approaches to
  turbulent multiphase flows},
\newblock \bibinfo{journal}{International Journal of Multiphase Flow}
  \bibinfo{volume}{35} (\bibinfo{year}{2009}) \bibinfo{pages}{801--810}.
%Type = Article
\bibitem[{Joseph et~al.(1990)Joseph, Lundgren, Jackson, and
  Saville}]{joseph1990}
\bibinfo{author}{D.~D. Joseph}, \bibinfo{author}{T.~S. Lundgren},
  \bibinfo{author}{R.~Jackson}, \bibinfo{author}{D.~A. Saville},
\newblock \bibinfo{title}{Ensemble averaged and mixture theory equations for
  incompressible fluid—particle suspensions},
\newblock \bibinfo{journal}{International journal of multiphase flow}
  \bibinfo{volume}{16} (\bibinfo{year}{1990}) \bibinfo{pages}{35--42}.
%Type = Book
\bibitem[{Simpson(1997)}]{simpson1997}
\bibinfo{author}{J.~E. Simpson}, \bibinfo{title}{Gravity Currents. In the
  environment and in the laboratory}, \bibinfo{publisher}{Cambridge University
  Press}, \bibinfo{edition}{2\textsuperscript{nd}} edition,
  \bibinfo{year}{1997}.
%Type = Article
\bibitem[{Zang(1991)}]{Zang1991}
\bibinfo{author}{T.~A. Zang},
\newblock \bibinfo{title}{{On the rotation and skew-symmetric forms for
  incompressible flow simulations}},
\newblock \bibinfo{journal}{Applied Numerical Mathematics} \bibinfo{volume}{7}
  (\bibinfo{year}{1991}) \bibinfo{pages}{27--40}.
%Type = Article
\bibitem[{Morinishi et~al.(1998)Morinishi, Lund, Vasilyev, and
  Moin}]{Morinishi1998}
\bibinfo{author}{Y.~Morinishi}, \bibinfo{author}{T.~Lund},
  \bibinfo{author}{O.~Vasilyev}, \bibinfo{author}{P.~Moin},
\newblock \bibinfo{title}{{Fully conservative higher order finite difference
  schemes for incompressible flow}},
\newblock \bibinfo{journal}{Journal of Computational Physics}
  \bibinfo{volume}{143} (\bibinfo{year}{1998}) \bibinfo{pages}{90--124}.
%Type = Article
\bibitem[{R{\o}nquist(1996)}]{Ronquist1996}
\bibinfo{author}{E.~M. R{\o}nquist},
\newblock \bibinfo{title}{{Convection treatment using spectral elements of
  different order}},
\newblock \bibinfo{journal}{International Journal for Numerical Methods in
  Fluids} \bibinfo{volume}{22} (\bibinfo{year}{1996})
  \bibinfo{pages}{241--264}.
%Type = Article
\bibitem[{Gatski(1991)}]{Gatski1991}
\bibinfo{author}{T.~B. Gatski},
\newblock \bibinfo{title}{Review of incompressible fluid flow computations
  using the vorticity-velocity formulation},
\newblock \bibinfo{journal}{Applied Numerical Mathematics} \bibinfo{volume}{7}
  (\bibinfo{year}{1991}) \bibinfo{pages}{227--239}.
%Type = Article
\bibitem[{Daube(1992)}]{Daube1992}
\bibinfo{author}{O.~Daube},
\newblock \bibinfo{title}{{Resolution of the 2D Navier-Stokes equations in
  velocity-vorticity form by means of an influence matrix technique}},
\newblock \bibinfo{journal}{Journal of Computational Physics}
  \bibinfo{volume}{103} (\bibinfo{year}{1992}) \bibinfo{pages}{402--414}.
%Type = Article
\bibitem[{Clercx(1997)}]{Clercx1997a}
\bibinfo{author}{H.~Clercx},
\newblock \bibinfo{title}{{A spectral solver for the Navier-Stokes equations in
  the velocity-vorticity formulation for flows with two nonperiodic
  directions}},
\newblock \bibinfo{journal}{Journal of Computational Physics}
  \bibinfo{volume}{137} (\bibinfo{year}{1997}) \bibinfo{pages}{186--211}.
%Type = Article
\bibitem[{Arnold et~al.(2010)Arnold, Falk, and Winther}]{arnold2010finite}
\bibinfo{author}{D.~N. Arnold}, \bibinfo{author}{R.~S. Falk},
  \bibinfo{author}{R.~Winther},
\newblock \bibinfo{title}{{Finite element exterior calculus: from Hodge theory
  to numerical stability}},
\newblock \bibinfo{journal}{Bulletin of the American Mathematical Society}
  \bibinfo{volume}{47} (\bibinfo{year}{2010}) \bibinfo{pages}{281--354}.
%Type = Article
\bibitem[{Palha et~al.(2014)Palha, Rebelo, Hiemstra, Kreeft, and
  Gerritsma}]{Palha2014}
\bibinfo{author}{A.~Palha}, \bibinfo{author}{P.~P. Rebelo},
  \bibinfo{author}{R.~Hiemstra}, \bibinfo{author}{J.~Kreeft},
  \bibinfo{author}{M.~I. Gerritsma},
\newblock \bibinfo{title}{Physics-compatible discretization techniques on
  single and dual grids, with application to the {Poisson} equation of volume
  forms},
\newblock \bibinfo{journal}{Journal of Computational Physics}
  \bibinfo{volume}{257} (\bibinfo{year}{2014}) \bibinfo{pages}{1394--1422}.
%Type = Article
\bibitem[{Bossavit(1999{\natexlab{a}})}]{bossavit_japan_computational_1}
\bibinfo{author}{A.~Bossavit},
\newblock \bibinfo{title}{{Computational electromagnetism and geometry: (1)
  Network equations}},
\newblock \bibinfo{journal}{Journal of the Japan Society of Applied
  Electromagnetics} \bibinfo{volume}{7} (\bibinfo{year}{1999}{\natexlab{a}})
  \bibinfo{pages}{150--159}.
%Type = Article
\bibitem[{Bossavit(1999{\natexlab{b}})}]{bossavit_japan_computational_2}
\bibinfo{author}{A.~Bossavit},
\newblock \bibinfo{title}{{Computational electromagnetism and geometry: (2)
  Network constitutive laws}},
\newblock \bibinfo{journal}{Journal of the Japan Society of Applied
  Electromagnetics} \bibinfo{volume}{7} (\bibinfo{year}{1999}{\natexlab{b}})
  \bibinfo{pages}{294--301}.
%Type = Article
\bibitem[{Bossavit(1999{\natexlab{c}})}]{bossavit_japan_computational_3}
\bibinfo{author}{A.~Bossavit},
\newblock \bibinfo{title}{{Computational electromagnetism and geometry: (3)
  Convergence}},
\newblock \bibinfo{journal}{Journal of the Japan Society of Applied
  Electromagnetics} \bibinfo{volume}{7} (\bibinfo{year}{1999}{\natexlab{c}})
  \bibinfo{pages}{401--408}.
%Type = Article
\bibitem[{Bossavit(2000{\natexlab{a}})}]{bossavit_japan_computational_4}
\bibinfo{author}{A.~Bossavit},
\newblock \bibinfo{title}{{Computational electromagnetism and geometry: (4)
  From degrees of freedom to fields}},
\newblock \bibinfo{journal}{Journal of the Japan Society of Applied
  Electromagnetics} \bibinfo{volume}{8} (\bibinfo{year}{2000}{\natexlab{a}})
  \bibinfo{pages}{102--109}.
%Type = Article
\bibitem[{Bossavit(2000{\natexlab{b}})}]{bossavit_japan_computational_5}
\bibinfo{author}{A.~Bossavit},
\newblock \bibinfo{title}{{Computational electromagnetism and geometry: (5) The
  ``Galerkin Hodge''}},
\newblock \bibinfo{journal}{Journal of the Japan Society of Applied
  Electromagnetics} \bibinfo{volume}{8} (\bibinfo{year}{2000}{\natexlab{b}})
  \bibinfo{pages}{203--209}.
%Type = Incollection
\bibitem[{Kirby et~al.(2012)Kirby, Logg, Rognes, and Terrel}]{kirby2012}
\bibinfo{author}{C.~Kirby}, \bibinfo{author}{A.~Logg}, \bibinfo{author}{M.~E.
  Rognes}, \bibinfo{author}{A.~R. Terrel},
\newblock \bibinfo{title}{Common and unusual finite elements},
\newblock in: \bibinfo{booktitle}{Automated Solution of Differential Equations
  by the Finite Element Method}, volume~\bibinfo{volume}{84} of
  \textit{\bibinfo{series}{Lecture Notes in Computational Science and
  Engineering}}, \bibinfo{publisher}{Springer}, \bibinfo{year}{2012}, pp.
  \bibinfo{pages}{95--119}.
%Type = Incollection
\bibitem[{Raviart and Thomas(2006)}]{RaviartThomas1977}
\bibinfo{author}{P.~Raviart}, \bibinfo{author}{J.~M. Thomas},
\newblock \bibinfo{title}{{A mixed finite element method for 2nd order elliptic
  problems}},
\newblock in: \bibinfo{booktitle}{Mathematical Aspects of the Finite Element
  Method}, volume \bibinfo{volume}{606} of \textit{\bibinfo{series}{Lecture
  Notes in Mathematics}}, \bibinfo{publisher}{Springer Berlin Heidelberg},
  \bibinfo{year}{2006}, pp. \bibinfo{pages}{292--315}.
%Type = Book
\bibitem[{Hairer et~al.(2006)Hairer, Lubich, and Wanner}]{Hairer2006}
\bibinfo{author}{E.~Hairer}, \bibinfo{author}{C.~Lubich},
  \bibinfo{author}{G.~Wanner}, \bibinfo{title}{Geometric numerical integration:
  structure-preserving algorithms for ordinary differential equations},
  volume~\bibinfo{volume}{31}, \bibinfo{publisher}{Springer Science \& Business
  Media}, \bibinfo{year}{2006}.
%Type = Book
\bibitem[{Brezzi and Fortin(1991)}]{brezzi1991}
\bibinfo{author}{F.~Brezzi}, \bibinfo{author}{M.~Fortin}, \bibinfo{title}{Mixed
  and Hybrid Finite Element Methods}, number~\bibinfo{number}{15} in
  \bibinfo{series}{Springer Series in Computational Mathematics 44},
  \bibinfo{publisher}{Springer}, \bibinfo{year}{1991}.

\end{thebibliography}




%	\bibitem{}

% 	\end{thebibliography}

\end{document}